\documentclass[11pt]{article} 
\usepackage{amsfonts,amsmath,latexsym,amssymb,mathrsfs,amsthm,comment}
\usepackage{graphicx}
\usepackage{xcolor}
\usepackage{booktabs}
\DeclareRobustCommand{\stirling}{\genfrac\{\}{0pt}{}}

\let\OLDthebibliography\thebibliography
\renewcommand\thebibliography[1]{
  \OLDthebibliography{#1}
  \setlength{\parskip}{1pt}
  \setlength{\itemsep}{0pt plus 0.0ex}
}

\def\numberlikeadb{\global\def\theequation{\thesection.\arabic{equation}}}
\numberlikeadb
\newtheorem{theorem}{Theorem}[section]
\newtheorem{lemma}[theorem]{Lemma}
\newtheorem{corollary}[theorem]{Corollary}

\newtheorem{proposition}[theorem]{Proposition}
\newtheorem{remark}[theorem]{Remark}

\usepackage{lscape}
\usepackage{caption}
\usepackage{multirow}
\allowdisplaybreaks
\begin{document}

\title{Asymptotic approximations for the distribution of the product of correlated normal random variables
%with application to risk measures
%Asymptotic approximations of  tail probabilities and quantile functions for the product of correlated normal random variables
%The distribution of the product of correlated normal random variables: tail behaviour and related distributional properties
}
\author{Robert E. Gaunt\footnote{Department of Mathematics, The University of Manchester, Oxford Road, Manchester M13 9PL, UK, robert.gaunt@manchester.ac.uk; zixin.ye@postgrad.manchester.ac.uk}\:\, and Zixin Y$\mathrm{e}^{*}$}

\date{} 
\maketitle

\vspace{-5mm}

\begin{abstract}  We obtain asymptotic approximations for the probability density function of the product of two correlated normal random variables with non-zero means and arbitrary variances. As a consequence, we deduce asymptotic approximations for the tail probabilities and quantile functions of this distribution, as well as an asymptotic approximation for the widely used risk measures value at risk and tail value at risk. 
\end{abstract}



\section{Introduction}

Let $(X, Y)$ be a bivariate normal random vector with mean vector $(\mu_X,\mu_Y)$, variances $(\sigma_X^2,\sigma_Y^2)$ and correlation coefficient $\rho$. The distribution of the product $Z=XY$ arises in numerous applications, with recent examples including condensed matter physics \cite{ach}, astrophysics \cite{cac} and chemical physics \cite{hey}. Since the work of \cite{craig,wb32} in the 1930's, the distributional theory of the product $Z=XY$ has also received much attention, 
%for example, formulas for the moment generating function, first four moments and skewness and kurtosis are given by\cite{craig}, the cumulant generating function and cumulants of general order are given by \cite{h42}
%with a number of the most basic and important distributional properties worked out in references such as
%with contributions including the works of 
with contributions coming from, amongst others, \cite{a47,bc08,h42}; see \cite{gaunt22,np16} for an overview of the literature, as well as \cite{gaunt22} for a review of the basic distributional theory in the zero mean case $\mu_X=\mu_Y=0$.
%Since the work of \cite{craig} in 1936, the distribution of the product $Z=XY$ has received much attention in the statistics literature (see \cite{gaunt22,np16} for an overview of some of the literature), and has found numerous applications, with recent examples including condensed matter physics \cite{ach}, astrophysics \cite{cac} and chemical physics \cite{hey}. 

Despite the interest in the distribution of the product $Z$, it was not until 2016 
%-- some 80 years after the work of \cite{craig} -- 
that \cite{cui} succeeded in deriving the following exact formula for the probability density function (PDF) of $Z$ that is valid for general $\mu_X,\mu_Y\in\mathbb{R}$, $\sigma_X,\sigma_Y>0$ and $-1<\rho<1$. For $x\in\mathbb{R}$,
\begin{align}f(x)&=\frac{1}{\pi}\exp\bigg\{-\frac{1}{2(1-\rho^2)}\bigg(\frac{\mu_X^2}{\sigma_X^2}+\frac{\mu_Y^2}{\sigma_Y^2}-\frac{2\rho(x+\mu_X\mu_Y)}{\sigma_X\sigma_Y}\bigg)\bigg\}\nonumber\\
&\quad\times\sum_{n=0}^\infty\sum_{m=0}^{2n}\frac{x^{2n-m}|x|^{m-n}}{(2n)!(\sigma_X\sigma_Y)^{n+1}(1-\rho^2)^{2n+1/2}}\binom{2n}{m}\bigg(\frac{\mu_X}{\sigma_X}-\frac{\rho \mu_Y}{\sigma_Y}\bigg)^m\nonumber\\
\label{pdf}&\quad\times\bigg(\frac{\mu_Y}{\sigma_Y}-\frac{\rho \mu_X}{\sigma_X}\bigg)^{2n-m}K_{m-n}\bigg(\frac{|x|}{(1-\rho^2)\sigma_X\sigma_Y}\bigg),
\end{align}
where $K_\nu(x)$ is a modified Bessel function of the second kind (see Appendix \ref{appa} for a definition). If $\rho=0$ and one of the means is equal to zero, then the PDF (\ref{pdf}) simplifies to a single infinite series (see \cite{cui,simon}). Without loss of generality take $\mu_Y=0$. Then the PDF (\ref{pdf}) simplifies to 
\begin{equation}\label{asdh}
f(x)=\frac{1}{\pi\sigma_X\sigma_Y}\exp\bigg(\!-\frac{\mu_X^2}{2\sigma_X^2}\bigg)\sum_{n=0}^\infty \frac{\mu_X^{2n}|x|^n}{(2n)!\sigma_X^{3n}\sigma_Y^n} K_n\bigg(\frac{|x|}{\sigma_X\sigma_Y}\bigg),\quad x\in\mathbb{R}.   
\end{equation}
It should be noted that formula (\ref{asdh}) is a Neumann–$K$ series of modified Bessel functions
of the second kind, for which there are methods for their closed-form summation; see, for example, \cite{BJMP}. Also, in the case $\mu_X=\mu_Y=0$, the following simpler exact formula for the PDF of $Z$ is available:
\begin{equation}f(x)=\frac{1}{\pi\sigma_X\sigma_Y\sqrt{1-\rho^2}}\exp\bigg(\frac{\rho  x}{\sigma_X\sigma_Y(1-\rho^2)} \bigg)K_0\bigg(\frac{ |x|}{\sigma_X\sigma_Y(1-\rho^2)}\bigg), \quad x\in\mathbb{R}, \label{pdf0}
\end{equation}
which has been derived independently by \cite{g96,np16,wb32}, and identified as a variance-gamma distribution by \cite{gaunt prod}. 
%The PDF in the case $\mu_X=\mu_Y=\rho=0$ was in fact known much earlier to \cite{craig}, as well as \cite{wb32,y33}. 
%We also note that the PDF (\ref{pdf}) simplifies to a single infinite series if either $\mu_X=0$ or $\mu_Y=0$ (see \cite{cui,simon}), but that there is not a significant simplification in the case that $\rho=0$ and $\mu_X$ and $\mu_Y$ are non-zero.
% and it was identified by \cite{gaunt prod} that $\overline{Z}_n$ is variance-gamma distributed. 
%However, an exact formula for the PDF of the mean $\overline{Z}_n$ is not known for the case of a general mean vector $(\mu_X,\mu_Y)\in\mathbb{R}^2$.

The exact PDF (\ref{pdf}) takes a rather complicated form, and so it is of particular interest to study the asymptotic properties of the distribution of the product $Z$. This serves as motivation for this paper, in which we study in detail the asymptotic properties of this distribution. We begin by obtaining asymptotic approximations for the PDF (\ref{pdf}) in the limits $x\rightarrow0$ and $|x|\rightarrow\infty$ (Proposition \ref{prop2} and Theorem \ref{prop3}). From the asymptotic behaviour as $x\rightarrow0$, we deduce that the distribution of the product $Z$ is unimodal with mode at 0 for all parameter constellations.  From our asymptotic approximations for the PDF in the limit $|x|\rightarrow\infty$, we deduce asymptotic approximations of the tail probabilities of the product $Z$ (Theorem \ref{corsec3}). Given the complicated form of the PDF (\ref{pdf}), these asymptotic approximations may be useful in allowing researchers to obtain approximate tail probabilities in a computationally efficient and simple manner. These asymptotic approximations generalise approximations of \cite{bc08,gaunt22} for the distribution of the product of two correlated zero mean normal random variables, and complement asymptotic approximations for the PDF (see \cite{g17,springer}) and tail probabilities (see \cite{l23}) of the product of $n$ independent zero mean normal random variables.

From our asymptotic approximations for the tail probabilities, we obtain asymptotic approximations for the quantile function of the product $Z$ (Theorem \ref{thmq}); note that a closed-form formula for the quantile function is not available. Quantile functions are important in statistics, for example in assessing statistical significance and in Monte Carlo simulations, and in mathematical finance in which the widely used risk measures value at risk (VaR) and tail value at risk (TVaR) are expressed in terms of the quantile function (see \cite{a99,mc15}). Indeed, our asymptotic approximations for the quantile function of the product $Z$ immediately yield asymptotic approximations for VaR, and we deduce asymptotic approximations for TVaR (Corollary \ref{cortvar}). 

These asymptotic approximations for VaR and TVaR, together with the tail probability estimates of Theorem \ref{corsec3}, could be of interest in portfolio analysis when considering the compound return of (possibly correlated) normally distributed one-period returns. Indeed, consider a portfolio with initial value $v_0$ at time 0, and suppose that the rate of returns $R_1$ and $R_2$ in years 1 and 2 are $N(\mu_1,\sigma_1^2)$
and $N(\mu_2,\sigma_2^2)$ distributed, respectively, with correlation $\rho$. Then the value of the portfolio $V$ at the end of the second year is given by $V=v_0(1+R_1)(1+R_2)=_d v_0 X_1X_2$, where $X_1\sim N(1+\mu_1,\sigma_1^2)$ and $X_2\sim (1+\mu_2,\sigma_2^2)$ are normal random variables with correlation $\rho$.
%We remark that in portfolio analysis, one-period returns are often assumed to be normally distributed, but that this assumption is not valid for multi-periods due to the effect of compounding.
%, in which it is more appropriate to use the distribution of the pro.

%(note that asset prices are often correlated, and that this feature is particularly important when assessing the risk of heavy losses in a portfolio).

%These asymptotic approximations could be of interest in portfolio analysis when considering the compound return of two assets with correlated normally distributed one-period returns (note that asset prices are often correlated, and that this feature is particularly important when assessing the risk of heavy losses in a portfolio).

The rest of the paper is organised as follows. In Section \ref{sec2.1}, we state our asymptotic approximation for the distribution of the product $Z$, whilst in Section \ref{sec2.2} we state our asymptotic approximations for VaR and TVaR. In Section \ref{sec2.3}, we present numerical results in order to assess the quality of the asymptotic approximations. %In %Section \ref{sec2} we state our main results. 
In Section \ref{sec3}, we state and prove some lemmas that are used in Section \ref{sec4}, in which we prove the main results from Section \ref{sec2}. Lastly, Appendix \ref{appa} collects some basic results on special functions that are used in this paper.

\section{Main results}\label{sec2}

\subsection{Asymptotic approximations for the distribution}\label{sec2.1}

The following Proposition \ref{prop2} and Theorem \ref{prop3} provide asymptotic approximations for the PDF (\ref{pdf}) of the product $Z=XY$. To allow for more compact formulas, we let
\begin{align*}C=C(\mu_X,\mu_Y,\sigma_X,\sigma_Y,\rho)=\exp\bigg\{-\frac{1}{2(1-\rho^2)}\bigg(\frac{\mu_X^2}{\sigma_X^2}+\frac{\mu_Y^2}{\sigma_Y^2}-\frac{2\rho\mu_X\mu_Y}{\sigma_X\sigma_Y}\bigg)\bigg\}.
\end{align*}

\begin{comment}
\begin{theorem}\label{thm1}
\begin{align*}F(x)=
\end{align*}
\end{theorem}
\end{comment}

\begin{proposition}\label{prop2}Let $\mu_X,\mu_Y\in\mathbb{R}$, $\sigma_X,\sigma_Y>0$ and $-1<\rho<1$. Then, as $x\rightarrow0$,
\begin{align}\label{xzero}f(x)\sim -\frac{C\ln|x|}{\pi\sigma_X\sigma_Y\sqrt{1-\rho^2}}.
\end{align}
In particular, the distribution of $Z$ is unimodal with mode $0$ for all parameter constellations.
\end{proposition}

\begin{theorem}\label{prop3} Let $\mu_X,\mu_Y\in\mathbb{R}$, $\sigma_X,\sigma_Y>0$ and $-1<\rho<1$. Then, as $x\rightarrow\infty$,
\begin{align}f(x)&= \frac{C}{\sqrt{2\pi\sigma_X\sigma_Y x}}\exp\bigg\{\frac{1}{8}\bigg(\frac{1+\rho}{1-\rho}\bigg)\bigg(\frac{\mu_X}{\sigma_X}-\frac{\mu_Y}{\sigma_Y}\bigg)^2-\frac{x}{\sigma_X\sigma_Y(1+\rho)}\bigg\}\nonumber\\
\label{for1}&\quad\times\cosh\bigg(\bigg|\frac{\mu_X}{\sigma_X}+\frac{\mu_Y}{\sigma_Y}\bigg|\frac{\sqrt{x}}{(1+\rho)\sqrt{\sigma_X\sigma_Y}}\bigg)\big\{1+O(x^{-1/2})\big\},
\end{align}
and, as $x\rightarrow-\infty$,
\begin{align}f(x)&= \frac{C}{\sqrt{2\pi\sigma_X\sigma_Y |x|}}\exp\bigg\{\frac{1}{8}\bigg(\frac{1-\rho}{1+\rho}\bigg)\bigg(\frac{\mu_X}{\sigma_X}+\frac{\mu_Y}{\sigma_Y}\bigg)^2+\frac{x}{\sigma_X\sigma_Y(1-\rho)}\bigg\}\nonumber\\
\label{for2}&\quad\times\cosh\bigg(\bigg|\frac{\mu_X}{\sigma_X}-\frac{\mu_Y}{\sigma_Y}\bigg|\frac{\sqrt{|x|}}{(1-\rho)\sqrt{\sigma_X\sigma_Y}}\bigg)\big\{1+O(|x|^{-1/2})\big\}.
\end{align}
If $\mu_X/\sigma_X+\mu_Y/\sigma_Y=0$, then the error in the asymptotic approximation (\ref{for1}) is of the smaller order $O(x^{-1})$. If $\mu_X/\sigma_X-\mu_Y/\sigma_Y=0$, then the error in the approximation (\ref{for2}) is of the smaller order $O(x^{-1})$. 
\end{theorem}

From Theorem \ref{prop3}, we deduce the following asymptotic approximations for the tail probabilities of the product $Z$. Let $F(x)=\mathbb{P}(Z\leq x)$ and $\bar{F}(x)=1-F(x)=\mathbb{P}(Z>x)$. In the following theorem, we find it more natural to state our asymptotic approximation for the tail probability of the product $Z$ in the limit $x\rightarrow\infty$ in terms of the survival function $\bar{F}(x)$ than in terms of the cumulative distribution function (CDF). For example, in deriving our asymptotic approximation for the quantile function of the product $Z$ in the limit $p\rightarrow1$ given in Theorem \ref{thmq}, we will apply formula (\ref{for3}) for the asymptotic approximation of the survival probability in the limit $x\rightarrow\infty$.
%We will later apply these formulas to derive asymptotic approximations for the quantile function of the product $Z$.}

\begin{theorem}\label{corsec3}Let $\mu_X,\mu_Y\in\mathbb{R}$, $\sigma_X,\sigma_Y>0$ and $-1<\rho<1$. Then, as $x\rightarrow\infty$,
\begin{align}\bar{F}(x)&= \frac{C(1+\rho)\sqrt{\sigma_X\sigma_Y}}{\sqrt{2\pi x}}\exp\bigg\{\frac{1}{8}\bigg(\frac{1+\rho}{1-\rho}\bigg)\bigg(\frac{\mu_X}{\sigma_X}-\frac{\mu_Y}{\sigma_Y}\bigg)^2-\frac{x}{\sigma_X\sigma_Y(1+\rho)}\bigg\}\nonumber\\
\label{for3}&\quad\times\cosh\bigg(\bigg|\frac{\mu_X}{\sigma_X}+\frac{\mu_Y}{\sigma_Y}\bigg|\frac{\sqrt{x}}{(1+\rho)\sqrt{\sigma_X\sigma_Y}}\bigg)\big\{1+O(x^{-1/2})\big\},
\end{align}
and, as $x\rightarrow-\infty$,
\begin{align}F(x)&= \frac{C(1-\rho)\sqrt{\sigma_X\sigma_Y}}{\sqrt{2\pi|x|}}\exp\bigg\{\frac{1}{8}\bigg(\frac{1-\rho}{1+\rho}\bigg)\bigg(\frac{\mu_X}{\sigma_X}+\frac{\mu_Y}{\sigma_Y}\bigg)^2+\frac{x}{\sigma_X\sigma_Y(1-\rho)}\bigg\}\nonumber\\
\label{for4}&\quad\times\cosh\bigg(\bigg|\frac{\mu_X}{\sigma_X}-\frac{\mu_Y}{\sigma_Y}\bigg|\frac{\sqrt{|x|}}{(1-\rho)\sqrt{\sigma_X\sigma_Y}}\bigg)\big\{1+O(|x|^{-1/2})\big\}.
\end{align}
If $\mu_X/\sigma_X+\mu_Y/\sigma_Y=0$, then the error in the asymptotic approximation (\ref{for3}) is of the smaller order $O(x^{-1})$. If $\mu_X/\sigma_X-\mu_Y/\sigma_Y=0$, then the error in the approximation (\ref{for4}) is of the smaller order $O(x^{-1})$.
\end{theorem}

 For $0<p<1$, let $Q(p)=F^{-1}(p)$ denote the quantile function of the product $Z$. The following theorem provides asymptotic approximations for the quantile function. 

\begin{theorem}\label{thmq}Let $\mu_X,\mu_Y\in\mathbb{R}$, $\sigma_X,\sigma_Y>0$ and $-1<\rho<1$. Let $\delta$ be such that $\delta=0$ if $\mu_X/\sigma_X+\mu_Y/\sigma_Y=0$ and $\delta=1$ if $\mu_X/\sigma_X+\mu_Y/\sigma_Y\not=0$. Let $q=1-p$. Also, let
\begin{align}
G(p)&=G(\mu_X,\mu_Y,\sigma_X,\sigma_Y,\rho,p)\nonumber\\
&=\sigma_X\sigma_Y(1+\rho)\bigg\{\ln(1/q)+\frac{1}{\sqrt{1+\rho}}\bigg|\frac{\mu_X}{\sigma_X}+\frac{\mu_Y}{\sigma_Y}\bigg|\sqrt{\ln(1/q)}-\frac{1}{2}\ln(\ln(1/q))\nonumber\\
&\quad-\frac{1}{2}\ln\bigg(\frac{2\pi}{1+\rho}\bigg)+\frac{1}{4(1+\rho)}\bigg(\frac{\mu_X}{\sigma_X}+\frac{\mu_Y}{\sigma_Y}\bigg)^2+\frac{1}{8}\bigg(\frac{1+\rho}{1-\rho}\bigg)\bigg(\frac{\mu_X}{\sigma_X}-\frac{\mu_Y}{\sigma_Y}\bigg)^2\nonumber\\
\label{gfor}&\quad-\frac{1}{2(1-\rho^2)}\bigg(\frac{\mu_X^2}{\sigma_X^2}+\frac{\mu_Y^2}{\sigma_Y^2}-\frac{2\rho\mu_X\mu_Y}{\sigma_X\sigma_Y}\bigg)-\delta\ln(2)\bigg\}.
\end{align}
Then, as $p\rightarrow1$,
\begin{align}\label{q1}Q(p)&=G(p)+O\bigg(\frac{\ln(\ln(1/(1-p)))}{\sqrt{\ln(1/(1-p))}}\bigg).
\end{align}
Now let $\delta'$ be such that $\delta'=0$ if $\mu_X/\sigma_X-\mu_Y/\sigma_Y=0$ and $\delta'=1$ if $\mu_X/\sigma_X-\mu_Y/\sigma_Y\not=0$. Then, as $p\rightarrow0$,
\begin{align}\label{q2}Q(p)&=-\sigma_X\sigma_Y(1-\rho)\bigg\{\ln(1/p)+\frac{1}{\sqrt{1-\rho}}\bigg|\frac{\mu_X}{\sigma_X}-\frac{\mu_Y}{\sigma_Y}\bigg|\sqrt{\ln(1/p)}-\frac{1}{2}\ln(\ln(1/p))\nonumber\\
&\quad-\frac{1}{2}\ln\bigg(\frac{2\pi}{1-\rho}\bigg)+\frac{1}{4(1-\rho)}\bigg(\frac{\mu_X}{\sigma_X}-\frac{\mu_Y}{\sigma_Y}\bigg)^2+\frac{1}{8}\bigg(\frac{1-\rho}{1+\rho}\bigg)\bigg(\frac{\mu_X}{\sigma_X}+\frac{\mu_Y}{\sigma_Y}\bigg)^2\nonumber\\
&\quad-\frac{1}{2(1-\rho^2)}\bigg(\frac{\mu_X^2}{\sigma_X^2}+\frac{\mu_Y^2}{\sigma_Y^2}-\frac{2\rho\mu_X\mu_Y}{\sigma_X\sigma_Y}\bigg)-\delta'\ln(2)\bigg\}+O\bigg(\frac{\ln(\ln(1/p))}{\sqrt{\ln(1/p)}}\bigg).
\end{align}
For $\delta=0$, the error in the asymptotic approximation (\ref{q1}) is of the smaller order $O(1/\ln(1/(1-p)))$. If $\delta'=0$, then the error in the approximation (\ref{q2}) is of the smaller order $O(1/\ln(1/p))$.
\end{theorem}

\begin{remark}
\noindent 1. From Proposition \ref{prop2} we observe that the density of the product $Z$ is locally symmetric about the mode at the origin. 

\vspace{2mm}

\noindent 2. Observe that the asymptotic approximation (\ref{for2}) can be obtained from (\ref{for1}) by replacing $(\mu_Y,\rho,x)$ by $(-\mu_Y,-\rho,-x)$. The asymptotic approximation (\ref{for4}) can be obtained from (\ref{for3}) using the same procedure. As noted in the proof of Theorem \ref{prop3}, this follows from the fact that $Z=XY=_d-X'Y'$, where $(X', Y')$ is a bivariate normal random vector with mean vector $(\mu_X,-\mu_Y)$, variances $(\sigma_X^2,\sigma_Y^2)$ and correlation coefficient $-\rho$. 

\vspace{2mm}

\noindent 3. We now discuss the asymptotic approximation (\ref{for1}); similar comments apply to the approximations (\ref{for2})--(\ref{for4}). There is a significant change in the asymptotic behaviour of the PDF as $x\rightarrow\infty$ depending on whether the quantity $\Delta=\mu_X/\sigma_X+\mu_Y/\sigma_Y$ is zero or non-zero. Indeed, if $\Delta=0$, then $f(x)\sim Ax^{-1/2}\mathrm{e}^{-ax}$, as $x\rightarrow\infty$, for some constants $A,a>0$, whilst if $\Delta\not=0$, then $f(x)\sim Ax^{-1/2}\mathrm{e}^{-ax+b\sqrt{x}}$, as $x\rightarrow\infty$, for some constants $A,a,b>0$. In particular, the tails become heavier as $|\Delta|$ increases.

We express the asymptotic approximation (\ref{for1}) in terms of the hyperbolic cosine function $\cosh(x)$ in order to give a compact formula that is valid for all possible values of $\Delta$. If one wishes to give an asymptotic approximation purely in terms of power functions and the exponential function, then the cases $\Delta=0$ and $\Delta\not=0$ would need to be treated separately. This is because 
%when $\delta=0$, the hyperbolic cosine factor disappears because $\cosh(0)=1$, whilst 
for $\Delta\not=0$ we have that $\cosh(\Delta\sqrt{x})=(\mathrm{e}^{\Delta\sqrt{x}}+\mathrm{e}^{-\Delta\sqrt{x}})/2=\mathrm{e}^{\Delta\sqrt{x}}/2+O(x^{-1/2})$ as $x\rightarrow\infty$, meaning we pick up a multiplicative factor of $1/2$ rather than a multiplicative factor of $\cosh(0)=1$ in the $\Delta=0$ case.

\vspace{2mm}

\noindent 4. The asymptotic approximations (\ref{q1}) and (\ref{q2}) contain the first four terms in asymptotic expansion of the quantile function $Q(p)$ in the limits $p\rightarrow1$ and $p\rightarrow0$, respectively. This is in spite of the fact that the approximations are derived as a consequence of the approximations of Theorem \ref{corsec3} in which just the leading term of the asymptotic expansion is present. 
%The asymptotic approximations (\ref{q1}) and (\ref{q2}) do, however, converge very slowly due to the presence of an error term that has a logarithmic-type order. 

We also observe that the leading term in the asymptotic expansions (\ref{q1}) and (\ref{q2}) does not involve the mean parameters $(\mu_X,\mu_Y)$; these parameters first appear in the second term in the expansion.
\end{remark}

\subsection{Application to risk measures}\label{sec2.2}

In mathematical finance, risk measures are used to determine the amount of assets to be kept in reserve. Risk measures are often used by regulators to ensure financial institutions have reasonable capital buffers to reduce risks of heavy losses and defaults. For further details on the applications and mathematical properties of risk measures see the seminal paper \cite{a99} and the influential book \cite{mc15}. The widely used risk measures value at risk (VaR) and tail value at risk (TVaR) (also known as the expected shortfall) are expressed in terms of quantile functions. As such, we can apply the asymptotic approximations from Theorem \ref{thmq} to derive asymptotic approximations for VaR and TVaR for the product $Z$, which could arise, for example, as the compound return of (possibly correlated) normally distributed one-period returns. 

Let us first recall the definitions of these risk measures.
Let $P$ be a random variable denoting the price change of an asset over a specified time interval. We define the \emph{loss} random variable to be $L=-P$. We will subsequently be considering the product $Z$ as a loss random variable.

The VaR at level $p \in (0, 1)$ of a loss random variable $L$
with CDF $F_L$ is defined as
\[\mathrm{VaR}_p(L):=F_L^{-1}(p).\]
It is the minimum loss in the $100(1-p)\%$ of worst case scenarios.

The TVaR at level $p \in (0, 1)$ of a loss random variable $L$ is defined by
\[\mathrm{TVaR}_p(L):=\mathbb{E}[L\,|\,L\geq\mathrm{VaR}_p(L)].\]
It is the expected loss in the $100(1-p)\%$ of worst case scenarios. 
%Note that for any level $p\in(0,1)$ and any loss random variable $L$, $\mathrm{TVaR}_p(L)\geq\mathrm{VaR}_p(L)$. 
If $L$ is a continuous random variable, then we have the useful formula
\[\mathrm{TVaR}_p(L)=\frac{1}{1-p}\int_p^1\mathrm{VaR}_t(L)\,\mathrm{d}t.\]
Here we have followed the convention used in actuarial applications, by defining VaR and TVaR in terms of a loss random variable, although we note that VaR and TVaR are also alternatively defined in terms of the profit random variable $P$ (see \cite{a99}).
%by $\mathrm{VaR}_p^{\mathrm{left}}(P)=F^{-1}_P(1-p)$ and $\mathrm{TVaR}_p^{\mathrm{left}}(P)=\mathbb{E}[-P\,|\,P\leq -\mathrm{VaR}(P)]$.

\begin{comment}
Let us first recall the definitions of these risk measures.
Let $P$ be a random variable denoting the price change of an asset over a specified time interval. We define the \emph{loss} random variable to be $L=-P$.

The VaR at level $p \in (0, 1)$ of a loss random variable $L$
with CDF $F_L$ is defined as
\[\mathrm{VaR}_p(L):=F_L^{-1}(p).\]
It is the minimum loss in the $100(1-p)\%$ of worst case scenario.

The TVaR at level $p \in (0, 1)$ of a loss random variable $L$ is defined by
\[\mathrm{TVaR}_p(L)=\mathbb{E}[L\,|\,L\geq\mathrm{VaR}_p(L)].\]
It is the expected loss in the $100(1-p)\%$ of worst case scenario. If $L$ is a continuous random variable, then we have the useful formula
\[\mathrm{TVaR}_p(L)=\frac{1}{1-p}\int_p^1\mathrm{VaR}_t(L)\,\mathrm{d}t.\]
Here we have followed the convention followed in actuarial applications, by defining VaR and TVaR in terms of a loss random variable, although we note that VaR and TVaR are also alternatively defined in terms of the profit random variable $P$ by $\mathrm{VaR}_p^{\mathrm{left}}(P)=F^{-1}_P(1-p)$ and $\mathrm{TVaR}_p^{\mathrm{left}}(P)=\mathbb{E}[-P\,|\,P\leq -\mathrm{VaR}(P)]$.
\end{comment}

\begin{corollary}\label{cortvar}Let $\mu_X,\mu_Y\in\mathbb{R}$, $\sigma_X,\sigma_Y>0$ and $-1<\rho<1$. Then, as $p\rightarrow1$,
\begin{align}\label{varfor}\mathrm{VaR}_p(Z)&=G(p)+O\bigg(\frac{\ln(\ln(1/(1-p)))}{\sqrt{\ln(1/(1-p))}}\bigg),\\
\label{tvarfor}\mathrm{TVaR}_p(Z)&=G(p)+\sigma_X\sigma_Y(1+\rho)+O\bigg(\frac{\ln(\ln(1/(1-p)))}{\sqrt{\ln(1/(1-p))}}\bigg),
\end{align}
where $G(p)=G(\mu_X,\mu_Y,\sigma_X,\sigma_Y,\rho,p)$ is defined as in (\ref{gfor}). For $\mu_X\sigma_X+\mu_Y/\sigma_Y=0$, the error in the asymptotic approximation (\ref{q1}) is of the smaller order $O(1/\ln(1/(1-p)))$.
%\begin{align}\label{tvarfor}
%\mathrm{TVaR}_p(Z)&= -\sigma_X\sigma_Y(1-\rho)\bigg\{\ln(1/p)+\frac{1}{\sqrt{1-\rho}}\bigg|\frac{\mu_X}{\sigma_X}-\frac{\mu_Y}{\sigma_Y}\bigg|\sqrt{\ln(1/p)}-\frac{1}{2}\ln(\ln(1/p))\nonumber\\
%&\quad+1-\frac{1}{2}\ln\bigg(\frac{2\pi}{1-\rho}\bigg)+\frac{1}{4}\bigg(\frac{\mu_X}{\sigma_X}-\frac{\mu_Y}{\sigma_Y}\bigg)^2+\frac{1}{8}\bigg(\frac{1-\rho}{1+\rho}\bigg)^2\bigg(\frac{\mu_X}{\sigma_X}+\frac{\mu_Y}{\sigma_Y}\bigg)^2\nonumber\\
%&\quad-\frac{1}{2(1-\rho^2)}\bigg(\frac{\mu_X^2}{\sigma_X^2}+\frac{\mu_Y^2}{\sigma_Y^2}-\frac{2\rho\mu_X\mu_Y}{\sigma_X\sigma_Y}\bigg)\bigg\}+o(1).
%\end{align}
\end{corollary}

\begin{remark}
From (\ref{varfor}) and (\ref{tvarfor}) we deduce that, as $p\rightarrow1$,
\begin{equation*}
\mathrm{TVaR}_p(Z)-\mathrm{VaR}_p(Z)=\sigma_X\sigma_Y(1+\rho)+o(1),
\end{equation*}
and
\begin{equation*}
\frac{\mathrm{TVaR}_p(Z)}{\mathrm{VaR}_p(Z)}-1=\frac{1}{\ln(1/(1-p))}\big(1+o(1)\big),  
\end{equation*}
and so the relative error between $\mathrm{VaR}_p(Z)$ and $\mathrm{TVaR}_p(Z)$ tends to zero as $p\rightarrow1$. We remark that this property is shared by the exponential distribution, for which exact formulas are available for VaR and TVaR. Indeed, for $L\sim \mathrm{Exp}(\lambda)$ with PDF $f_L(x)=\lambda\mathrm{e}^{-\lambda x}$, $x>0$, we have $\mathrm{VaR}_p(L)=\lambda^{-1}\ln(1/(1-p))$, $\mathrm{TVaR}_p(L)=(\ln(1/(1-p))+1)/\lambda$ (see \cite{varef}), and so $\mathrm{TVaR}_p(L)/\mathrm{VaR}_p(L)-1=1/\ln(1/(1-p))$.
\end{remark}

\subsection{Numerical results}\label{sec2.3}

%Suppose $\sigma_X=\sigma_Y=1$. Let us express the product $Z=XY$ in terms of independent standard normal random variables. We first note that we can write $Z=_d(U+\mu_X)(V+\mu_Y)$, where $U$ and $V$ are independent standard normal random variables with correlation $\rho$. Let $W=(V-\rho U)/\sqrt{1-\rho^2}$, then it is readily seen that $U$ and $W$ are independent standard normal random variables. We therefore have that $Z=_d(U+\mu_X)(\rho U+\sqrt{1-\rho^2}W+\mu_Y)$, where $U$ and $W$ are independent standard normal random variables.

\begin{table}[h]
  \centering
  \caption{\footnotesize{Relative error in approximating the PDF of $Z$ by the asymptotic approximation (\ref{for1}). Here $\sigma_X=\sigma_Y=1$.
  A negative number means that the approximation
is less than the true value.
  }}
\label{table1}
\footnotesize{
\begin{tabular}{l*{6}{c}}
\hline
& \multicolumn{6}{c}{$x$} \\
\cmidrule(lr){2-7}
$(\mu_X,\mu_Y,\rho)$ & 2.5 & 5 & 7.5 & 10 & 12.5 & 15 \\
\hline
        (0,0,-0.5) & 3.4E-02 & 1.8E-02 & 1.2E-02 & 9.1E-03 & 7.3E-03 & 6.1E-03 \\ 
        (0,0,-0.25) & 4.1E-02 & 2.2E-02 & 1.5E-02 & 1.1E-02 & 9.1E-03 & 7.6E-03 \\ 
        (0,0,0) & 4.4E-02 & 2.3E-02 & 1.6E-02 & 1.2E-02 & 9.7E-03 & 8.1E-03 \\ 
        (0,0,0.25) & 4.1E-02 & 2.2E-02 & 1.5E-02 & 1.1E-02 & 9.1E-03 & 7.6E-03 \\ 
        (0,0,0.5) & 3.4E-02 & 1.8E-02 & 1.2E-02 & 9.1E-03 & 7.3E-03 & 6.1E-03 \\ \hline
        (1,-1,-0.5) & 4.4E-02 & 2.3E-02 & 1.6E-02 & 1.2E-02 & 9.7E-03 & 8.1E-03 \\
        (1,-1,-0.25) & 6.3E-02 & 3.4E-02 & 2.3E-02 & 1.8E-02 & 1.4E-02 & 1.2E-02 \\ 
        (1,-1,0) & 8.2E-02 & 4.5E-02 & 3.1E-02 & 2.4E-02 & 1.9E-02 & 1.6E-02 \\
        (1,-1,0.25) & 1.0E-01 & 5.5E-02 & 3.8E-02 & 2.9E-02 & 2.4E-02 & 2.0E-02 \\
        (1,-1,0.5) & 1.2E-01 & 6.7E-02 & 4.6E-02 & 3.5E-02 & 2.9E-02 & 2.4E-02 \\ \hline
        (2,-2,-0.5) & 7.4E-02 & 4.0E-02 & 2.7E-02 & 2.1E-02 & 1.7E-02 & 1.4E-02 \\
        (2,-2,-0.25) & 1.3E-01 & 7.0E-02 & 4.9E-02 & 3.7E-02 & 3.0E-02 & 2.5E-02 \\ 
        (2,-2,0) & 1.9E-01 & 1.1E-01 & 7.5E-02 & 5.7E-02 & 4.7E-02 & 3.9E-02 \\ 
        (2,-2,0.25) & 2.7E-01 & 1.5E-01 & 1.1E-01 & 8.2E-02 & 6.7E-02 & 5.6E-02 \\
        (2,-2,0.5) & 3.7E-01 & 2.1E-01 & 1.4E-01 & 1.1E-01 & 9.0E-02 & 7.6E-02 \\ \hline
        (1,0,-0.5) & -8.8E-02 & -7.1E-02 & -6.1E-02 & -5.5E-02 & -5.0E-02 & -4.6E-02 \\ 
        (1,0,-0.25) & -5.8E-02 & -5.3E-02 & -4.8E-02 & -4.4E-02 & -4.0E-02 & -3.8E-02 \\ 
        (1,0,0) & -3.2E-02 & -3.8E-02 & -3.6E-02 & -3.4E-02 & -3.2E-02 & -3.0E-02 \\
        (1,0,0.25) & -1.1E-02 & -2.4E-02 & -2.5E-02 & -2.5E-02 & -2.4E-02 & -2.3E-02 \\ 
        (1,0,0.5) & 4.9E-03 & -1.3E-02 & -1.6E-02 & -1.7E-02 & -1.7E-02 & -1.7E-02 \\ \hline
        (1,1,-0.5) & -2.1E-01 & -1.6E-01 & -1.3E-01 & -1.1E-01 & -1.0E-01 & -9.4E-02 \\ 
        (1,1,-0.25) & -1.5E-01 & -1.2E-01 & -1.0E-01 & -9.0E-02 & -8.1E-02 & -7.5E-02 \\ 
        (1,1,0) & -1.1E-01 & -8.8E-02 & -7.6E-02 & -6.8E-02 & -6.2E-02 & -5.7E-02 \\ 
        (1,1,0.25) & -7.1E-02 & -6.1E-02 & -5.3E-02 & -4.8E-02 & -4.4E-02 & -4.1E-02 \\ 
        (1,1,0.5) & -4.0E-02 & -3.7E-02 & -3.3E-02 & -3.0E-02 & -2.8E-02 & -2.6E-02 \\ \hline
        (2,1,-0.5) & -3.9E-01 & -2.8E-01 & -2.3E-01 & -2.0E-01 & -1.8E-01 & -1.6E-01 \\ 
        (2,1,-0.25) & -3.1E-01 & -2.3E-01 & -1.9E-01 & -1.7E-01 & -1.5E-01 & -1.4E-01 \\ 
        (2,1,0) & -2.4E-01 & -1.9E-01 & -1.6E-01 & -1.4E-01 & -1.2E-01 & -1.1E-01 \\ 
        (2,1,0.25) & -1.9E-01 & -1.5E-01 & -1.3E-01 & -1.1E-01 & -1.0E-01 & -9.4E-02 \\ 
        (2,1,0.5) & -1.4E-01 & -1.1E-01 & -9.9E-02 & -8.8E-02 & -8.0E-02 & -7.5E-02 \\ \hline
    \end{tabular}}
\end{table}

\begin{table}[h]
  \centering
  \caption{\footnotesize{Relative error in approximating the tail probability $\mathbb{P}(Z\geq x)$ by the asymptotic approximation (\ref{for3}). Here $\sigma_X=\sigma_Y=1$.
  A negative number means that the approximation
is less than the true value. N/A means that the asymptotic approximation is not valid, in which case no result is reported.
  }}
\label{table2}
\footnotesize{
\begin{tabular}{l*{6}{c}}
\hline
& \multicolumn{6}{c}{$x$} \\
\cmidrule(lr){2-7}
$(\mu_X,\mu_Y,\rho)$ & $q_{0.95}$ & $q_{0.975}$ & $q_{0.99}$ & $q_{0.995}$ & $q_{0.999}$ & $q_{0.9999}$ \\
\hline
        (0, 0, -0.5) & 3.6E-01 & 2.7E-01 & 2.0E-01 & 1.7E-01 & 1.2E-01 & 8.8E-02 \\ 
        (0, 0, -0.25) & 3.1E-01 & 2.4E-01 & 1.8E-01 & 1.6E-01 & 1.1E-01 & 8.1E-02 \\ 
        (0, 0, 0) & 2.8E-01 & 2.2E-01 & 1.7E-01 & 1.4E-01 & 1.1E-01 & 7.7E-02 \\ 
        (0, 0, 0.25) & 2.6E-01 & 2.0E-01 & 1.6E-01 & 1.3E-01 & 1.0E-01 & 7.3E-02 \\ 
        (0, 0, 0.5) & 2.3E-01 & 1.8E-01 & 1.4E-01 & 1.2E-01 & 9.3E-02 & 6.8E-02 \\ \hline
        (1, -1, -0.5) & 5.2E-01 & 3.6E-01 & 2.5E-01 & 2.1E-01 & 1.5E-01 & 1.0E-01 \\ 
        (1, -1, -0.25) & 4.6E-01 & 3.3E-01 & 2.4E-01 & 2.0E-01 & 1.4E-01 & 9.9E-02 \\ 
        (1, -1, 0) & 4.3E-01 & 3.1E-01 & 2.3E-01 & 1.9E-01 & 1.4E-01 & 9.7E-02 \\ 
        (1, -1, 0.25) & 4.1E-01 & 3.0E-01 & 2.2E-01 & 1.9E-01 & 1.4E-01 & 9.7E-02 \\ 
        (1, -1, 0.5) & 4.0E-01 & 3.0E-01 & 2.2E-01 & 1.9E-01 & 1.4E-01 & 9.8E-02 \\ \hline
        (2, -2, -0.5) & N/A & 1.5E+00 & 6.0E-01 & 4.2E-01 & 2.5E-01 & 1.6E-01 \\ 
        (2, -2, -0.25) & N/A & 1.3E+00 & 6.0E-01 & 4.3E-01 & 2.6E-01 & 1.7E-01 \\ 
        (2, -2, 0) & N/A & 1.3E+00 & 6.2E-01 & 4.5E-01 & 2.8E-01 & 1.8E-01 \\ 
        (2, -2, 0.25) & N/A & 1.3E+00 & 6.5E-01 & 4.8E-01 & 2.9E-01 & 1.9E-01 \\ 
        (2, -2, 0.5) & N/A & 1.4E+00 & 6.9E-01 & 5.1E-01 & 3.2E-01 & 2.1E-01 \\ \hline
        (1, 0, -0.5) & -3.0E-01 & -2.9E-01 & -2.7E-01 & -2.6E-01 & -2.4E-01 & -2.2E-01 \\ 
        (1, 0, -0.25) & -2.0E-01 & -2.0E-01 & -2.0E-01 & -2.0E-01 & -1.8E-01 & -1.7E-01 \\ 
        (1, 0, 0) & -1.4E-01 & -1.5E-01 & -1.5E-01 & -1.5E-01 & -1.5E-01 & -1.4E-01 \\ 
        (1, 0, 0.25) & -8.8E-02 & -1.0E-01 & -1.1E-01 & -1.2E-01 & -1.2E-01 & -1.1E-01 \\ 
        (1, 0, 0.5) & -5.3E-02 & -7.4E-02 & -8.6E-02 & -9.1E-02 & -9.4E-02 & -9.2E-02 \\ \hline
        (1, 1, -0.5) & -5.5E-01 & -5.2E-01 & -4.9E-01 & -4.7E-01 & -4.4E-01 & -4.0E-01 \\ 
        (1, 1, -0.25) & -4.6E-01 & -4.3E-01 & -4.1E-01 & -3.9E-01 & -3.6E-01 & -3.3E-01 \\ 
        (1, 1, 0) & -3.9E-01 & -3.7E-01 & -3.5E-01 & -3.3E-01 & -3.0E-01 & -2.8E-01 \\ 
        (1, 1, 0.25) & -3.3E-01 & -3.2E-01 & -3.0E-01 & -2.9E-01 & -2.6E-01 & -2.4E-01 \\ 
        (1, 1, 0.5) & -2.9E-01 & -2.7E-01 & -2.6E-01 & -2.5E-01 & -2.3E-01 & -2.1E-01 \\ \hline
        (2, 1, -0.5) & -6.9E-01 & -6.6E-01 & -6.3E-01 & -6.1E-01 & -5.7E-01 & -5.3E-01 \\ 
        (2, 1, -0.25) & -6.1E-01 & -5.8E-01 & -5.5E-01 & -5.3E-01 & -4.9E-01 & -4.5E-01 \\ 
        (2, 1, 0) & -5.5E-01 & -5.2E-01 & -4.9E-01 & -4.7E-01 & -4.3E-01 & -4.0E-01 \\ 
        (2, 1, 0.25) & -4.9E-01 & -4.7E-01 & -4.4E-01 & -4.2E-01 & -3.9E-01 & -3.6E-01 \\ 
        (2, 1, 0.5) & -4.5E-01 & -4.2E-01 & -4.0E-01 & -3.8E-01 & -3.5E-01 & -3.2E-01 \\ \hline
    \end{tabular}}
\end{table}

\begin{table}[h]
  \centering
  \caption{\footnotesize{Relative error in approximating the quantile function $Q(p)$ by the asymptotic approximation (\ref{q1}).
  Here $\sigma_X=\sigma_Y=1$.
  A negative number means that the approximation
is less than the true value.
N/A means that the asymptotic approximation is not valid, in which case no result is reported.  }}
\label{table3}
\footnotesize{
\begin{tabular}{l*{6}{c}}
\hline
& \multicolumn{6}{c}{$p$} \\
\cmidrule(lr){2-7}
$(\mu_X,\mu_Y,\rho)$ &        0.95 &    0.975 & 0.99 &  0.995 &  0.999 & 0.9999 \\
\hline
(0, 0, -0.5)	&	-1.7E-01	&	-8.6E-02	&	-5.4E-02	&	-3.4E-02	&	-2.5E-02	&	-1.6E-02	\\
(0, 0, -0.25)	&	-9.9E-02	&	-5.7E-02	&	-3.8E-02	&	-2.5E-02	&	-1.9E-02	&	-1.2E-02	\\
(0, 0, 0)	&	-6.5E-02	&	-4.2E-02	&	-3.0E-02	&	-2.0E-02	&	-1.7E-02	&	-1.0E-02	\\
(0, 0, 0.25)	&	-4.8E-02	&	-3.3E-02	&	-2.5E-02	&	-1.8E-02	&	-1.4E-02	&	-9.8E-03	\\
(0, 0, 0.5)	&	-3.7E-02	&	-2.8E-02	&	-2.2E-02	&	-1.6E-02	&	-1.3E-02	&	-9.1E-03	\\
\hline
(1, -1, -0.5)	&	-6.7E-01	&	-2.2E-01	&	-1.2E-01	&	-6.3E-02	&	-4.4E-02	&	-2.2E-02	\\
(1, -1, -0.25)	&	-3.6E-01	&	-1.5E-01	&	-8.2E-02	&	-4.7E-02	&	-3.5E-02	&	-2.0E-02	\\
(1, -1, 0)	&	-2.4E-01	&	-1.1E-01	&	-6.3E-02	&	-3.8E-02	&	-2.8E-02	&	-1.6E-02	\\
(1, -1, 0.25)	&	-1.6E-01	&	-8.0E-02	&	-4.9E-02	&	-3.1E-02	&	-2.4E-02	&	-1.3E-02	\\
(1, -1, 0.5)	&	-1.1E-01	&	-6.1E-02	&	-3.9E-02	&	-2.5E-02	&	-2.0E-02	&	-1.2E-02	\\
\hline
(2, -2, -0.5)	&	N/A	&	3.4E+00	&	-2.0E+00	&	-3.7E-01	&	-1.9E-01	&	-7.3E-02	\\
(2, -2, -0.25)	&	N/A	&	5.2E+00	&	-1.1E+00	&	-2.6E-01	&	-1.4E-01	&	-5.8E-02	\\
(2, -2, 0)	&	N/A	&	5.7E+00	&	-7.1E-01	&	-1.9E-01	&	-1.1E-01	&	-4.7E-02	\\
(2, -2, 0.25)	&	N/A	&	4.5E+00	&	-4.9E-01	&	-1.4E-01	&	-8.4E-02	&	-3.5E-02	\\
(2, -2, 0.5)	&	N/A	&	2.8E+00	&	-3.3E-01	&	-1.0E-01	&	-6.1E-02	&	-2.8E-02	\\
\hline
(1, 0, -0.5)	&	1.6E-01	&	9.2E-02	&	6.0E-02	&	3.8E-02	&	2.7E-02	&	1.5E-02	\\
(1, 0, -0.25)	&	1.2E-01	&	7.6E-02	&	5.3E-02	&	3.6E-02	&	2.9E-02	&	1.7E-02	\\
(1, 0, 0)	&	8.8E-02	&	6.1E-02	&	4.5E-02	&	3.1E-02	&	2.6E-02	&	1.6E-02	\\
(1, 0, 0.25)	&	6.1E-02	&	4.6E-02	&	3.6E-02	&	2.7E-02	&	2.2E-02	&	1.5E-02	\\
(1, 0, 0.5)	&	3.8E-02	&	3.3E-02	&	2.8E-02	&	2.3E-02	&	1.9E-02	&	1.3E-02	\\
\hline
(1, 1, -0.5)	&	-5.6E-02	&	-6.4E-02	&	-6.6E-02	&	-6.5E-02	&	-6.3E-02	&	-5.8E-02	\\
(1, 1, -0.25)	&	-7.8E-03	&	-2.2E-02	&	-2.9E-02	&	-3.3E-02	&	-3.3E-02	&	-3.4E-02	\\
(1, 1, 0)	&	1.9E-02	&	4.3E-04	&	-8.5E-03	&	-1.5E-02	&	-1.7E-02	&	-1.9E-02	\\
(1, 1, 0.25)	&	3.6E-02	&	1.5E-02	&	4.2E-03	&	-3.2E-03	&	-6.4E-03	&	-1.0E-02	\\
(1, 1, 0.5)	&	4.7E-02	&	2.5E-02	&	1.3E-02	&	4.6E-03	&	5.3E-04	&	-4.1E-03	\\
\hline
(2, 1, -0.5)	&	-2.1E-01	&	-2.0E-01	&	-1.8E-01	&	-1.7E-01	&	-1.6E-01	&	-1.4E-01	\\
(2, 1, -0.25)	&	-1.4E-01	&	-1.3E-01	&	-1.3E-01	&	-1.2E-01	&	-1.1E-01	&	-9.8E-02	\\
(2, 1, 0)	&	-9.8E-02	&	-9.5E-02	&	-9.1E-02	&	-8.5E-02	&	-8.1E-02	&	-7.3E-02	\\
(2, 1, 0.25)	&	-6.6E-02	&	-6.8E-02	&	-6.7E-02	&	-6.4E-02	&	-6.1E-02	&	-5.6E-02	\\
(2, 1, 0.5)	&	-4.2E-02	&	-4.8E-02	&	-4.9E-02	&	-4.8E-02	&	-4.7E-02	&	-4.3E-02	\\
\hline
    \end{tabular}}
\end{table}

\begin{table}[h]
  \centering
  \caption{\footnotesize{Relative error in approximating TVaR by the asymptotic approximation (\ref{tvarfor}).
  Here $\sigma_X=\sigma_Y=1$.
  A negative number means that the approximation
is less than the true value.
N/A means that the asymptotic approximation is not valid, in which case no result is reported.   }}
\label{table4}
\footnotesize{
\begin{tabular}{l*{6}{c}}
\hline
& \multicolumn{6}{c}{$p$} \\
\cmidrule(lr){2-7}
$(\mu_X,\mu_Y,\rho)$ &        0.95 &    0.975 & 0.99 &  0.995 &  0.999 & 0.9999 \\
\hline

(0, 0, -0.5)	&	3.3E-02	&	1.5E-02	&	7.3E-03	&	3.2E-03	&	1.6E-03	&	-3.5E-04	\\
(0, 0, -0.25)	&	4.6E-02	&	2.4E-02	&	1.3E-02	&	6.4E-03	&	4.2E-03	&	2.0E-03	\\
(0, 0, 0)	&	5.2E-02	&	2.7E-02	&	1.6E-02	&	8.2E-03	&	5.9E-03	&	2.1E-03	\\
(0, 0, 0.25)	&	5.4E-02	&	2.9E-02	&	1.8E-02	&	1.0E-02	&	6.8E-03	&	2.8E-03	\\
(0, 0, 0.5)	&	5.5E-02	&	3.0E-02	&	1.8E-02	&	1.0E-02	&	6.5E-03	&	3.0E-03	\\
\hline
(1, -1, -0.5)	&	-2.6E-02	&	-1.8E-02	&	-1.5E-02	&	-9.9E-03	&	-7.9E-03	&	-5.6E-03	\\
(1, -1, -0.25)	&	9.9E-03	&	1.8E-03	&	-2.0E-03	&	-3.0E-03	&	-2.8E-03	&	-2.5E-03	\\
(1, -1, 0)	&	3.0E-02	&	1.3E-02	&	5.4E-03	&	2.2E-03	&	-9.8E-04	&	-1.7E-03	\\
(1, -1, 0.25)	&	4.3E-02	&	2.1E-02	&	1.1E-02	&	5.6E-03	&	2.9E-03	&	8.2E-05	\\
(1, -1, 0.5)	&	5.3E-02	&	2.7E-02	&	1.5E-02	&	8.1E-03	&	4.9E-03	&	2.9E-03	\\
\hline
(2, -2, -0.5)	&	N/A	 & -6.4E-01	&	-2.4E-01	&	-1.0E-01	&	-6.8E-02	&	-3.3E-02	\\
(2, -2, -0.25)	& N/A	&	-3.7E-01	&	-1.5E-01	&	-7.1E-02	&	-4.7E-02	&	-2.2E-02	\\
(2, -2, 0)	&	N/A	&	-2.2E-01	&	-9.1E-02	&	-4.5E-02	&	-3.1E-02	&	-1.5E-02	\\
(2, -2, 0.25)	&	N/A	&	-1.2E-01	&	-5.0E-02	&	-2.5E-02	&	-1.8E-02	&	-9.5E-03	\\
(2, -2, 0.5)	&	N/A	&	-2.8E-02	&	-1.4E-02	&	-9.3E-03	&	-6.9E-03	&	-5.1E-03	\\
\hline
(1, 0, -0.5)	&	-4.1E-03	&	-8.1E-03	&	-1.0E-02	&	-1.1E-02	&	-1.1E-02	&	-1.1E-02	\\
(1, 0, -0.25)	&	4.2E-03	&	4.2E-04	&	-1.3E-03	&	-3.3E-03	&	-3.9E-03	&	-4.4E-03	\\
(1, 0, 0)	&	6.3E-03	&	2.9E-03	&	1.7E-03	&	5.2E-04	&	-1.3E-03	&	-2.4E-03	\\
(1, 0, 0.25)	&	4.8E-03	&	3.1E-03	&	2.0E-03	&	1.1E-03	&	5.5E-04	&	5.5E-04	\\
(1, 0, 0.5)	&	2.2E-03	&	2.2E-03	&	1.6E-03	&	1.7E-03	&	1.0E-03	&	-8.8E-04	\\
\hline
(1, 1, -0.5)	&	-1.6E-01	&	-1.4E-01	&	-1.3E-01	&	-1.1E-01	&	-1.0E-01	&	-8.8E-02	\\
(1, 1, -0.25)	&	-1.2E-01	&	-1.0E-01	&	-9.0E-02	&	-7.9E-02	&	-7.2E-02	&	-6.1E-02	\\
(1, 1, 0)	&	-8.8E-02	&	-7.6E-02	&	-6.7E-02	&	-5.9E-02	&	-5.4E-02	&	-4.6E-02	\\
(1, 1, 0.25)	&	-6.7E-02	&	-5.9E-02	&	-5.2E-02	&	-4.6E-02	&	-4.2E-02	&	-3.5E-02	\\
(1, 1, 0.5)	&	-5.2E-02	&	-4.5E-02	&	-4.1E-02	&	-3.6E-02	&	-3.3E-02	&	-2.9E-02	\\
\hline
(2, 1, -0.5)	&	-2.9E-01	&	-2.6E-01	&	-2.3E-01	&	-2.1E-01	&	-1.9E-01	&	-1.7E-01	\\
(2, 1, -0.25)	&	-2.3E-01	&	-2.0E-01	&	-1.8E-01	&	-1.6E-01	&	-1.5E-01	&	-1.3E-01	\\
(2, 1, 0)	&	-1.9E-01	&	-1.6E-01	&	-1.5E-01	&	-1.3E-01	&	-1.2E-01	&	-1.0E-01	\\
(2, 1, 0.25)	&	-1.6E-01	&	-1.4E-01	&	-1.2E-01	&	-1.1E-01	&	-9.9E-02	&	-8.4E-02	\\
(2, 1, 0.5)	&	-1.4E-01	&	-1.2E-01	&	-1.0E-01	&	-9.2E-02	&	-8.4E-02	&	-7.0E-02	\\

\hline
    \end{tabular}}
\end{table}

Table \ref{table1} gives the relative error in approximating the PDF (\ref{pdf}) of the product $Z$ by the asymptotic approximation (\ref{for1}). Here and in all our numerical results we took $\sigma_X=\sigma_Y=1$. For given parameters and a given $x$, we used \emph{Mathematica} to compute the PDF (\ref{pdf}) by truncating the infinite series at $n=50$, which was accurate and fast to implement. Table \ref{table2} gives the relative error in approximating the tail probability $\mathbb{P}(Z\geq x)$ by the asymptotic approximation (\ref{for3}) when $x$ is taken to be a quantile $q_p$ for $p$ ranging from 0.95 to 0.9999. Tables \ref{table3} and \ref{table4} give the relative error in approximating the quantile function and TVaR of the product $Z$ by the asymptotic approximations (\ref{q1}) and (\ref{tvarfor}), respectively. We obtained the results in Tables \ref{table2}--\ref{table4} by performing Monte Carlo simulations using Python. We simulated the distribution of $Z$ by simulating two independent standard normal random variables and then used the relation $Z=_d(U+\mu_X)(\rho U+\sqrt{1-\rho^2}V+\mu_Y)$, where $U$ and $V$ are independent $N(0,1)$ random variables (in the case $\sigma_X=\sigma_Y=1$). In this manner, for given parameters,  we generated $10^{10}$ realisations of the distribution of $Z$ from which the results in each of Tables \ref{table2}--\ref{table4} were derived.
We took the $k$-th  largest realisation as the empirical quantile (where $k=\lfloor N(1-p)\rfloor +1$), which we used to estimate $Q(p)$. We estimated TVaR using the empirical TVaR, given by $k^{-1}\sum_{j=1}^k x_{j:N}$, where $x_{j:N}$ is the $(N-j+1)$-th largest observation from simulated values $x_{1},x_2,\ldots,x_N$ of $Z$. For some parameter regimes, the quantile $Q(0.95)$ is negative, in which case the asymptotic approximation (\ref{for3}) is not valid (it requires $x>0$). Since the asymptotic approximations (\ref{q1}) and (\ref{tvarfor}) for the quantile function and TVaR, respectively, are derived from the approximation (\ref{for3}) these approximations are also not valid for $p=0.95$. These instances are denoted by N/A in the Tables \ref{table2}--\ref{table4}.

It took 10 hours to produce Tables \ref{table2}--\ref{table4} using a computer server at The University of Manchester. Results are reported to 2 significant figures (s.f.) and most entries in the tables are indeed accurate to 2 s.f. (confirmed by repeating the simulations several times), although for $p=0.9999$ we found that some of the results in Table \ref{table4} for the case $\mu_X+\mu_Y=0$ were in some cases not accurate to 1 s.f.\ because the error from the simulations was of a non-negligible size compared to the asymptotic approximation error. Nevertheless, the presented results give a good account of the accuracy of the approximations across a range of parameter values and $p$-values.

From Table \ref{table1} we see that the relative error in approximating the PDF by the approximation (\ref{for1}) decreases roughly according to the $x^{-1/2}$ rate when $\mu_X+\mu_Y\not=0$, whilst it decreases roughly at a rate $x^{-1}$ when $\mu_X+\mu_Y=0$. This is consistent with the order of the error of the asymptotic approximation (\ref{for1}) as stated in Theorem \ref{prop3}. Similarly, from Tables \ref{table3} and \ref{table4} we that the asymptotic approximations (\ref{q1}) and (\ref{tvarfor}) for the quantile function and TVaR decrease at a faster rate with respect to $p$ when $\mu_X+\mu_Y=0$. 

We see from Table \ref{table1} that the relative error in approximating the PDF by the approximation (\ref{for1}) increases as the the magnitude of the means $\mu_X$ and $\mu_Y$ increase. We also see that for $(\mu_X,\mu_Y)=(1,-1),(2,-2)$ the relative error tends to increase as $\rho$ increases from $-0.5$ to $0.5$, whilst the reverse is true for $(\mu_X,\mu_Y)=(1,0),(1,1),(2,1)$.  Such trends do not seem to occur for the approximations of the quantile function and TVaR. We do, however, see that accuracy of the asymptotic approximations for the quantile function and TVaR are similar. For a given $p$-value, our approximations for the quantile function and TVaR tend to be more accurate than our approximation for the tail probability. This can be understood from the fact that while $p$ may be close to 1 (meaning the approximations (\ref{q1}) and (\ref{tvarfor}) are accurate) the quantitle $x=q_p$ may be not be large enough for the approximation (\ref{for3}) to be accurate. 
%We remark that for the parameter constellation $(\mu_X,\mu_Y,\rho)=(2,-2,0.5)$ our approximations for the quantile function and TVaR were poor for all $p$-values we considered with relative errors well in excess of 1. Otherwise, 
Overall, the accuracy of our asymptotic approximations for the PDF, tail probability, quantile function and TVaR are quite reasonable, especially given that they are very simple and based on just the leading order term in the asymptotic expansion for the PDF of $Z$.

%Overall, the accuracy of our asymptotic approximations for the PDF, tail probability, quantile function and TVaR were quite reasonable, except for the parameter constellations $(\mu_X,\mu_Y,\rho)=(2,-2,0.25),(2,-2,0.5)$, for which our approximations for the quantile function and TVaR were poor for all $p$-values we considered with relative errors well in excess of 1.

\section{Preliminary lemmas}\label{sec3}

%In this section, we state and prove some lemmas that are used in the proofs of our main results in Section \ref{sec4}. 
The following lemma collects some binomial identities from \cite{r96}.

\begin{lemma}\label{lem0} Let $n,k\geq0$ be non-negative integers. 
\begin{comment}
Let
\begin{equation*}S_{n,k}:=\sum_{m=0}^{2n}\binom{2n}{m}(-1)^m(m-n)^k.
\end{equation*}
Then
\begin{align*}S_{n,k}=\begin{cases} (2n)!, & \quad k=2n, \\
0, & \quad k<2n. \end{cases} 
\end{align*}
\end{comment}
Then
\begin{align*}
%S_{n,k}:=
\sum_{m=0}^{2n}\binom{2n}{m}(-1)^m(m-n)^k=\begin{cases} (2n)!, & \quad k=2n, \\
0, & \quad k<2n. \end{cases} 
\end{align*}
\end{lemma} 

\begin{lemma}\label{lem1} Let $ k,n\geq1$ be positive integers, and let $a,b\in\mathbb{R}$ be such that $a+b\not=0$. Then
\begin{equation*}\sum_{m=0}^{2n}\binom{2n}{m}a^mb^{2n-m}(m-n)^{2k}=(a+b)^{2n}\bigg[\bigg(\frac{a-b}{a+b}\bigg)^{2k}n^{2k}+p_{a,b,k}(n) \bigg],
\end{equation*}
where $p_{a,b,k}(n)$ is a polynomial in $n$ of degree $2k-1$ with coefficients involving the constants $a$, $b$ and $k$.
\begin{comment}
Let 
\begin{align*}S_{a,b,k,n}:=\frac{1}{(a+b)^{2n}}\sum_{m=0}^{2n}\binom{2n}{m}a^mb^{2n-m}(m-n)^{2k}.
\end{align*}
%Suppose first that $a\not=b$. 
Then, as $n\rightarrow\infty$,
\begin{align*}S_{a,b,k,n}=\bigg(\frac{a-b}{a+b}\bigg)^{2k}n^{2k} +O(n^{2k-1}).
\end{align*}
We emphasise that when $a=b$ we have $S_{a,a,k,n}=O(n^{2k-1})$, $n\rightarrow\infty$. 
%Now suppose $a=b$. Then, as $n\rightarrow\infty$,
%\begin{align*}S_{a,a,k,n}=O(n^{2k-1})
%\end{align*}
\end{comment}
\end{lemma}

In the following lemma, we recognise an infinite series as a generalized hypergeometric function; see Appendix \ref{appa} for a definition.

\begin{lemma}\label{lemhyp} Let $k\geq1$. Then, for $x\in\mathbb{R}$,
\begin{equation}\label{series}\sum_{n=0}^\infty\frac{n^{2k}}{(2n)!}x^n=\frac{x}{2}\,{}_{2k-1}F_{2k}\bigg({2,\ldots,2 \atop 1,\ldots,1,\frac{3}{2}};\frac{x}{4}\bigg).
\end{equation}
Moreover, as $x\rightarrow\infty$,
\begin{equation}\label{pfq}\sum_{n=0}^\infty\frac{n^{2k}}{(2n)!}x^n=\frac{1}{2^{2k+1}}x^{k}\mathrm{e}^{\sqrt{x}}\big(1+O(x^{-1/2})\big).
\end{equation}
\end{lemma}

\begin{lemma}\label{lem2} (i) Let $a>0$ and $b\in\mathbb{R}$. Then, as $x\rightarrow\infty$,
\begin{equation}\label{int1}\int_x^\infty\frac{1}{\sqrt{t}}\exp(-at+b\sqrt{t})\,\mathrm{d}t=\frac{1}{a\sqrt{x}}\exp(-ax+b\sqrt{x})\big(1+O(x^{-1/2})\big).
\end{equation}
If $b=0$, then the error in the asymptotic approximation (\ref{int1}) is of the smaller order $O(x^{-1})$.

\vspace{2mm}

\noindent (ii) Suppose that the functions $g_k$, $k=1,2$, are such that $g_k(x)=O(x^{-k/2})$ as $x\rightarrow\infty$. Then, as $x\rightarrow\infty$,
\[\int_x^\infty\frac{g_k(t)}{\sqrt{t}}\exp(-at+b\sqrt{t})\,\mathrm{d}t=O\bigg(\frac{1}{x^{1/2+k/2}}\exp(-ax+b\sqrt{x})\bigg), \quad k=1,2.\]
\end{lemma}

\begin{lemma}\label{lemq} Let $a,A,z>0$ and $b\in\mathbb{R}$. Let $g:(0,\infty)\rightarrow\mathbb{R}$ be a function such that $g(x)=O(x^{-1/2})$ as $x\rightarrow\infty$.  Consider the solution $x$ of the equation
\begin{equation}\label{zeqn}\frac{A}{\sqrt{x}}\exp(-ax+b\sqrt{x})\big(1+g(x)\big)=z,
\end{equation}
and note that for sufficiently small $z$ there is a unique solution. Then, as $z\rightarrow0$,
\begin{align}\label{lemqeqn}x=\frac{1}{a}\ln(1/z)+\frac{b}{a^{3/2}}\sqrt{\ln(1/z)}-\frac{1}{2a}\ln(\ln(1/z))+\frac{b^2}{4a^2}+\frac{\ln(A\sqrt{a})}{a}+O\bigg(\frac{\ln(\ln(1/z))}{\sqrt{\ln(1/z)}}\bigg).
\end{align}
\end{lemma}

\begin{remark}
In the case $b=0$ and $g(x)=0$ for all $x>0$, equation (\ref{zeqn}) can be solved in terms of the Lambert $W$ function, for which the asymptotic properties are well-understood; see \cite{c96}.   
\end{remark}

%\noindent{\emph{Proof of Lemma \ref{lem0}.}}

%\vspace{2mm}

\noindent{\emph{Proof of Lemma \ref{lem1}.}} Suppose $a+b\not=0$. Let 
\begin{align*}S_{a,b,k,n}:=\frac{1}{(a+b)^{2n}}\sum_{m=0}^{2n}\binom{2n}{m}a^mb^{2n-m}(m-n)^{2k}.
\end{align*}
 Let $p=a/(a+b)$. Then we can write
\begin{align}
S_{a,b,k,n}&=\sum_{m=0}^{2n}\binom{2n}{m}p^m(1-p)^{2n-m}(m-n)^{2k}\nonumber\\
\label{id0}&=\sum_{j=0}^{2k}\binom{2k}{j}(-n)^{2k-j}\sum_{m=0}^{2n}\binom{2n}{m}p^m(1-p)^{2m-n}m^j,
\end{align}
where we wrote out the binomial expansion of $(m-n)^{2k}$ and then changed the order of summation. By formula (2.9) of \cite{k08} we have that
\begin{equation}\label{id}
 \sum_{m=0}^{2n}\binom{2n}{m}p^m(1-p)^{2m-n}m^j=\sum_{i=0}^j  \stirling{j}{i}p^i(2n)^{\underline{i}},
\end{equation}
where $\stirling{j}{i}=\sum_{r=0}^i(-1)^{i-r}\binom{j}{r}r^j$ is a Stirling number of the second kind (see \cite{olver}) and $x^{\underline{d}}=x(x-1)\cdots(x-d+1)$ denotes the falling factorial. Note that the identity (\ref{id}) was derived in the context of formulas for moments of the binomial distribution in which $0<p<1$, although the argument used by \cite{k08} is valid for $p\in\mathbb{R}$. We can develop the right-hand side of (\ref{id}) to obtain that
\begin{align}
\label{id2}\sum_{m=0}^{2n}\binom{2n}{m}p^m(1-p)^{2m-n}m^j=\begin{cases} (2n)^jp^j+q_{p,j}(n), & j\geq1, \\
1, & j=0, \end{cases} 
\end{align}
where $q_{p,j}(n)$ is a polynomial in $n$ of degree $j-1$ with coefficients involving $p$ and $j$.
On substituting (\ref{id2}) into (\ref{id0}) we obtain that
\begin{align*}
S_{a,b,k,n}&=\sum_{j=0}^{2k}\binom{2k}{j}(-n)^{2k-j}  (2n)^jp^j+Q_{p,k}(n) =n^{2k}(2p-1)^{2k}+Q_{p,k}(n),
\end{align*}
where $Q_{p,k}(n)$ is a polynomial in $n$ of degree $2k-1$ with coefficients involving $p$ and $k$, and we evaluated the sum using the binomial theorem. The lemma now follows on writing the right-hand side of the above expression in terms of $a$ and $b$.
\begin{comment}
from which we observe that
\begin{equation*}
S_{a,b,k,n}=\mathbb{E}[(X-m)^{2k}],    
\end{equation*}
where $X$ follows the binomial distribution with parameters $2n$ and $p$ and probability mass function $\mathbb{P}(X=x)=\binom{2n}{x}p^x(1-p)^{2n-x}$, $x=0,1,\ldots,2n$.
\end{comment}
\hfill $\Box$

\vspace{2mm}

\noindent{\emph{Proof of Lemma \ref{lemhyp}.}} We begin by writing
\begin{equation}\label{dfgh}\sum_{n=0}^\infty\frac{n^{2k}}{(2n)!}x^n=\sum_{n=1}^\infty\frac{n^{2k}}{(2n)!}x^n=x\sum_{j=0}^\infty\frac{(j+1)^{2k}}{(2j+2)!}x^j.
\end{equation}
By the duplication formula \cite[Section 5.5(iii)]{olver}, we have 
\[(2j+2)!=\Gamma(2j+3)=\frac{1}{\sqrt{\pi}}2^{2j+2} j!\Gamma(j+3/2),\]
and substituting into (\ref{dfgh}) we get that
\begin{align*}\sum_{n=0}^\infty\frac{n^{2k}}{(2n)!}x^n=\frac{\sqrt{\pi}}{4}x\sum_{j=0}^\infty\frac{(j+1)^{2k}}{j!\Gamma(j+3/2)}\bigg(\frac{x}{4}\bigg)^j.
\end{align*}
Since $j+1=(j+1)!/j!=(2)_j/(1)_j$ and $\Gamma(j+3/2)=\Gamma(3/2)(3/2)_j=(\sqrt{\pi}/2)(3/2)_j$, equation (\ref{series}) follows from the series representation (\ref{gauss}) of the generalized hypergeometric function. Finally, the limiting form (\ref{pfq}) follows from applying the limiting form (\ref{gausslim}) to equation (\ref{series}). \hfill $\Box$

\vspace{2mm}

\noindent{\emph{Proof of Lemma \ref{lem2}.}}
% Using the definition (\ref{erfcdefn}) of the complementary error function $\mathrm{erfc}(x)$ as an integral
(i) A straightforward integration gives that
\begin{equation*}I_{a,b}(x):=\int_x^\infty\frac{1}{\sqrt{t}}\exp(-at+b\sqrt{t})\,\mathrm{d}t=\sqrt{\frac{\pi}{a}}\mathrm{e}^{b^2/(4a)}\mathrm{erfc}\bigg(\frac{2a\sqrt{x}-b}{2\sqrt{a}}\bigg),
\end{equation*}
where the complementary error function $\mathrm{erfc}(x)$ is defined in equation (\ref{erfcdefn}). Applying the asymptotic expansion (\ref{erfc}) now yields that, as $x\rightarrow\infty$,
\begin{align}I_{a,b}(x)&=\frac{\mathrm{e}^{b^2/(4a)}}{a\sqrt{x}-b/2}\exp\bigg(-\bigg(\frac{2a\sqrt{x}-b}{2\sqrt{a}}\bigg)^2\bigg)\big(1+O(x^{-1})\big) \nonumber\\
\label{int11}&=\frac{1}{a\sqrt{x}}\exp(-ax+b\sqrt{x})\big(1+O(x^{-1/2})\big).
\end{align}
Lastly, we note that it is clear that if $b=0$, then the error in the asymptotic approximation (\ref{int11}) is of the smaller order $O(x^{-1})$. 

\vspace{2mm}

\noindent(ii) This follows from a straightforward integration by parts. We omit the details. \hfill $\Box$

\vspace{2mm}

\noindent{\emph{Proof of Lemma \ref{lemq}.}} For ease of notation, we let $w=A/z$ and $h(x)=1+g(x)$. Rearranging equation (\ref{zeqn}) using basic properties of the natural logarithm gives that
\begin{align}\label{ax1}ax-b\sqrt{x}=\ln(w)-\frac{1}{2}\ln(x)+\ln(h(x)).
\end{align}
Applying the quadratic formula and taking the positive solution, we have that
\begin{align*}
\sqrt{x}=\frac{1}{2a}\Big(b+\sqrt{b^2+4a(\ln(w)-(1/2)\ln(x)+\ln(h(x))}\Big),
\end{align*}
so that
\begin{align}
x&=\frac{1}{4a^2}\Big(2b^2+2b\sqrt{b^2+4a(\ln(w)-(1/2)\ln(x)+\ln(h(x))}\nonumber\\
\label{ax2}&\quad+4a(\ln(w)-(1/2)\ln(x)+\ln(h(x))\Big).    
\end{align}
We now note that it is clear from equation (\ref{ax1}) that $x=O(\ln(w))$ as $w\rightarrow\infty$, from which we deduce that $\ln(x)=O(\ln(\ln(w)))$ as $w\rightarrow\infty$, and that $\ln(h(x))=O(1/\sqrt{\ln(w)})$ as $w\rightarrow\infty$ (recall that $\ln(1+u)=O(u)$ as $u\rightarrow0$). Applying the expansion $\sqrt{1-u}=1-u/2+O(u^2)$ as $u\rightarrow0$ 
%and $\ln(1+u)=1+O(u)$
to (\ref{ax2}) now gives that, as $w\rightarrow\infty$,
\begin{align}
x&=\frac{1}{4a^2}\bigg\{2b^2+\bigg[4b\sqrt{a}\sqrt{\ln(w)}-\frac{b}{2\sqrt{a}}\frac{2a\ln(x)-4a\ln(h(x))-b^2}{\sqrt{\ln(w)}}+O\bigg(\frac{(\ln(x))^2}{(\ln(w))^{3/2}}\bigg)\bigg]\nonumber\\
&\quad+4a\ln(w)-2a\ln(x)+O\bigg(\frac{1}{\sqrt{\ln(w)}}\bigg)\bigg\}
\nonumber\\
\label{ax3}&=\frac{1}{a}\ln(w)+\frac{b}{a^{3/2}}\sqrt{\ln(w)}+\frac{b^2}{2a^2}-\frac{1}{2a}\ln(x)+O\bigg(\frac{\ln(\ln(w))}{\sqrt{\ln(w)}}\bigg).
\end{align}
We now recursively apply the asymptotic approximation (\ref{ax3}) to obtain that, as $w\rightarrow\infty$,
\begin{align}\label{ax4}x=\frac{1}{a}\ln(w)+\frac{b}{a^{3/2}}\sqrt{\ln(w)}+\frac{b^2}{2a^2}-\frac{1}{2a}\ln\bigg(\frac{1}{a}\ln(w)+r(w)\bigg)+O\bigg(\frac{\ln(\ln(w))}{\sqrt{\ln(w)}}\bigg),
\end{align}
where $r(w)=O(\sqrt{\ln(w)})$ as $w\rightarrow\infty$. Using that $\ln(1+u)=O(u)$ as $u\rightarrow0$, we deduce that, as $w\rightarrow\infty$,
\begin{align}\label{ax5}
\ln\bigg(\frac{1}{a}\ln(w)+r(w)\bigg)&=\ln\bigg(\frac{1}{a}\ln(w)\bigg)+\ln\bigg(1+\frac{ar(w)}{\ln(w)}\bigg)\nonumber\\
&=\ln(\ln(w))-\ln(a)+O\bigg(\frac{1}{\sqrt{\ln(w)}}\bigg).
\end{align}
On inserting the asymptotic approximation (\ref{ax5}) into (\ref{ax4}), and recalling that $w=A/z$, we obtain that, as $z\rightarrow0$,
\begin{align}\label{lemqeqn2} x=\frac{1}{a}\ln(A/z)+\frac{b}{a^{3/2}}\sqrt{\ln(A/z)}-\frac{1}{2a}\ln(\ln(A/z))+\frac{b^2}{4a^2}+\frac{\ln(a)}{2a}+O\bigg(\frac{\ln(\ln(1/z))}{\sqrt{\ln(1/z)}}\bigg).
\end{align} 
Using our usual arguments, it is readily seen that, as $z\rightarrow0$, 
\begin{align*}
\sqrt{\ln(A/z)}=\sqrt{\ln(1/z)}\sqrt{1+\frac{\ln(A)}{\ln(1/z)}}=\sqrt{\ln(1/z)}+O\bigg(\frac{1}{\sqrt{\ln(1/z)}}\bigg),
\end{align*}
and
\begin{align*}
\ln(\ln(A/z))=\ln(\ln(1/z))+\ln\bigg(1+\frac{\ln(A)}{\ln(1/z)}\bigg)=\ln(\ln(1/z))+O\bigg(\frac{1}{\ln(1/z)}\bigg),
\end{align*}
and applying these asymptotic approximations to (\ref{lemqeqn2}) now yields the asymptotic approximation (\ref{lemqeqn}).
\hfill $\Box$

\section{Proofs of main results}\label{sec4}

To ease notation, we shall prove our results for the case $\sigma_X=\sigma_Y=1$, with the results for the general case $\sigma_X,\sigma_Y>0$ following from the basic fact that $Z=XY=_d \sigma_X\sigma_Y UV$, where $(U,V)$ is a bivariate normal random vector with mean vector $(\mu_X/\sigma_X,\mu_Y/\sigma_Y)$, variances $(1,1)$ and correlation coefficient $\rho$.

\begin{comment}
\vspace{2mm}

\noindent{\emph{Proof of Theorem \ref{thm1}.}}
\end{comment}

\vspace{2mm}

\noindent{\emph{Proof of Proposition \ref{prop2}.}} From the limiting form (\ref{Ktend0}), we see that the term in the double sum of (\ref{pdf}) with $m=n=0$ is of order $O(\ln|x|)$ as $x\rightarrow0$, whilst all other terms are of smaller order in this limit. Applying the limiting form (\ref{Ktend0}) to the term in the double sum with $m=n=0$ yields the limiting form (\ref{xzero}). From the limiting form (\ref{xzero}) we see that the PDF of $Z$ has a singularity at $x=0$, and since the modified Bessel function $K_\nu(x)$ is bounded for $x\not=0$, it follows that the density (\ref{pdf}) is bounded for $x\not=0$, and so has a unique mode at $x=0$ (where the density has a singularity). \hfill $\Box$

\vspace{2mm}

%\begin{align}\label{f11}f(x)=M_1(x)\sum_{n=0}^\infty\frac{(\mu_X(1+\rho))^{2n}x^n}{(2n)!(1-\rho^2)^{2n}}\sum_{m=0}^{2n}\binom{2n}{m}(-1)^mK_{m-n}\bigg(\frac{x}{1-\rho^2}\bigg),
%\end{align}
%where
%\begin{equation*}M_1(x)=\frac{1}{\pi\sqrt{1-\rho^2}}\exp\bigg\{-\frac{\mu_X^2}{1-\rho}+\frac{\rho x}{1-\rho^2}\bigg\}.
%\end{equation*}
% Applying the asymptotic expansion (\ref{Ktendinfinity}) to equation (\ref{f11}) gives that, as $x\rightarrow\infty$,

\noindent{\emph{Proof of Theorem \ref{prop3}.}} We first consider the limiting behaviour as $x\rightarrow\infty$.  Applying the asymptotic expansion (\ref{Ktendinfinity}) to the PDF (\ref{pdf}) gives that, as $x\rightarrow\infty$,
\begin{align}\label{first}f(x)&\sim M(x)\sum_{n=0}^\infty\sum_{m=0}^{2n}\frac{x^{n}}{(2n)!(1-\rho^2)^{2n}}\binom{2n}{m}(\mu_X-\rho \mu_Y)^m(\mu_Y-\rho \mu_X)^{2n-m}\nonumber\\
&\quad\times\sum_{k=0}^\infty a_k(m-n)\frac{(1-\rho^2)^k}{x^k},
\end{align}
where the constants $a_k(m-n)$, $k\geq0$, are defined as in (\ref{akv}), and
\begin{equation*}M(x)=\frac{1}{\sqrt{2\pi x}}\exp\bigg\{-\frac{\mu_X^2+\mu_Y^2-2\rho\mu_X\mu_Y}{2(1-\rho^2)}\bigg\}\exp\bigg(-\frac{x}{1+\rho}\bigg).
\end{equation*}

We shall now consider the cases  $\mu_X+\mu_Y=0$ and $\mu_X+\mu_Y\not=0$ separately. First, we suppose $\mu_X+\mu_Y=0$. Letting $\mu_Y=-\mu_X$ in (\ref{first}) gives that, as $x\rightarrow\infty$,
\begin{align*}f(x)\sim M(x)\sum_{n=0}^\infty\frac{(\mu_X(1+\rho))^{2n}x^n}{(2n)!(1-\rho^2)^{2n}}\sum_{m=0}^{2n}\binom{2n}{m}(-1)^m\sum_{k=0}^\infty a_k(m-n)\frac{(1-\rho^2)^k}{x^k}.
\end{align*}  
%where
%\begin{align*}M_2(x)=\frac{1}{\sqrt{2\pi x}}\exp\bigg\{-\frac{\mu_X^2}{1-\rho}-\frac{ x}{1+\rho}\bigg\}.
%\end{align*}
We now observe that $a_k(\nu)$ can be written as
\begin{align}\label{akm}a_k(\nu)=\frac{\nu^{2k}}{k!2^k}+b_k(\nu), \quad k\geq0,
\end{align}
where $b_k(\nu)=\sum_{j=0}^{k-1} A_j \nu^{2j}$, $k\geq0$,
for some constants $A_0,A_1,\ldots,A_{k-1}$. On applying Lemma \ref{lem0}, we now obtain that
\begin{align*}f(x)&=M(x)\sum_{n=0}^\infty\frac{(\mu_X(1+\rho))^{2n}x^n}{(2n)!(1-\rho^2)^{2n}}\sum_{m=0}^{2n}\binom{2n}{m}(-1)^m a_n(m-n)\frac{(1-\rho^2)^n}{x^n}\big\{1+O(x^{-1})\big\}\\
&=M(x)\sum_{n=0}^\infty\frac{\mu_X^{2n}}{(2n)!}\bigg(\frac{1+\rho}{1-\rho}\bigg)^n\sum_{m=0}^{2n}\binom{2n}{m}(-1)^m\frac{(m-n)^{2n}}{n!2^n}\big\{1+O(x^{-1})\big\}\\
&=M(x)\sum_{n=0}^\infty\frac{\mu_X^{2n}}{(2n)!}\bigg(\frac{1+\rho}{1-\rho}\bigg)^n\frac{(2n)!}{n!2^n}\big\{1+O(x^{-1})\big\}\\
&=M(x)\exp\bigg\{\frac{\mu_X^2}{2}\bigg(\frac{1+\rho}{1-\rho}\bigg)\bigg\}\big\{1+O(x^{-1})\big\},
\end{align*}
and so we have derived the asymptotic approximation (\ref{for1}) for the case $\mu_X+\mu_Y=0$.

Now suppose $\mu_X+\mu_Y\not=0$. Using equation (\ref{akm}), we write the asymptotic expansion (\ref{first}) in the following form: as $x\rightarrow\infty$,
\begin{align*}
f(x)&\sim f_1(x)+R_1(x),
\end{align*}
where
\begin{align*}f_1(x)&= M(x)\sum_{n=0}^\infty\sum_{k=0}^\infty\frac{(1-\rho^2)^{k-2n}x^{n-k}}{(2n)!k!2^k}\\
&\quad\times\sum_{m=0}^{2n}\binom{2n}{m}(\mu_X-\rho\mu_Y)^m(\mu_Y-\rho\mu_X)^{2n-m}(m-n)^{2k},\\
R_1(x)&=M(x)\sum_{n=0}^\infty\sum_{k=0}^\infty\frac{(1-\rho^2)^{k-2n}x^{n-k}}{(2n)!}\\
&\quad\times\sum_{m=0}^{2n}\binom{2n}{m}(\mu_X-\rho\mu_Y)^m(\mu_Y-\rho\mu_X)^{2n-m}b_k(m-n).
\end{align*}
We shall argue later that $R_1(x)/f_1(x)=O(x^{-1})$ as $x\rightarrow\infty$. Now we concentrate on the term $f_1(x)$. Applying Lemma \ref{lem1} (which is permitted since we have assumed that $\mu_X+\mu_Y\not=0$) gives that, as $x\rightarrow\infty$,
\begin{align*}f_1(x)\sim f_2(x)+R_2(x),
\end{align*} 
where
\begin{align}f_2(x)&= M(x)\sum_{n=0}^\infty\sum_{k=0}^\infty\frac{(1-\rho^2)^{k-2n}x^{n-k}}{(2n)!k!2^k} ((1-\rho)(\mu_X+\mu_Y))^{2n-2k}\nonumber\\
\label{ftwo}&\qquad\qquad\qquad\qquad\qquad\qquad\qquad\times((1+\rho)(\mu_X-\mu_Y))^{2k}n^{2k},\\
R_2(x)&=M(x)\sum_{n=0}^\infty\sum_{k=0}^\infty\frac{(1-\rho^2)^{k-2n}x^{n-k}}{(2n)!}((1-\rho)(\mu_X+\mu_Y))^{2n}c_{k}(n),\nonumber
\end{align}
where $c_0(n)=1$ and $c_{k}(n)$, $k\geq1$, $n\geq0$, are polynomials in $n$ of degree $2k-1$ with coefficients involving $k$, as well as $(\mu_X,\mu_Y,\rho)$, which we suppress in the notation.  We shall argue later that $R_2(x)/f_2(x)=O(x^{-1/2})$ as $x\rightarrow\infty$; we now focus on $f_2(x)$. 

Interchanging the order of summation in (\ref{ftwo}) gives that
%, as $x\rightarrow\infty$,
\begin{align*}f_2(x)&=M(x)\sum_{k=0}^\infty\frac{(1-\rho^2)^k}{k!(2x)^k}\bigg(\frac{(1+\rho)(\mu_X-\mu_Y)}{(1-\rho)(\mu_X+\mu_Y)}\bigg)^{2k}\sum_{n=0}^\infty\frac{n^{2k}}{(2n)!}\bigg(\frac{(\mu_X+\mu_Y)^2x}{(1+\rho)^2}\bigg)^n.
%\\
%&=\frac{M(x)}{2}\sum_{k=0}^\infty\frac{(1-\rho^2)^k}{k!(2x)^k}\bigg(\frac{(1+\rho)(\mu_X-\mu_Y)}{(1-\rho)(\mu_X+\mu_Y)}\bigg)^{2k}\\
%&\quad\times \frac{(\mu_X+\mu_Y)^2x}{(1+\rho)^2}\, {}_{2k-1}F_{2k}\bigg({2,\ldots,2 \atop 1,\ldots,1,\frac{3}{2}};\frac{(\mu_X+\mu_Y)^2x}{4(1+\rho)^2}\bigg),
\end{align*}
%where we used Lemma \ref{lemhyp} to put the infinite series over $n$ into the generalized hypergeometric series form (\ref{gauss}). 
Applying the asymptotic approximation (\ref{pfq}) and simplifying now gives that, as $x\rightarrow\infty$,
\begin{align}
f_2(x)&=\frac{M(x)}{2}\exp\bigg(\frac{(\mu_X+\mu_Y)\sqrt{x}}{1+\rho}\bigg)\sum_{k=0}^\infty   \frac{1}{k!}\bigg(\frac{1}{8}\bigg(\frac{1+\rho}{1-\rho}\bigg)(\mu_X-\mu_Y)^2\bigg)^k\big\{1+O(x^{-1/2})\big\}\nonumber\\
\label{f2}&=\frac{M(x)}{2}\exp\bigg\{\frac{1}{8}\bigg(\frac{1+\rho}{1-\rho}\bigg)(\mu_X-\mu_Y)^2\bigg\}\exp\bigg(\frac{(\mu_X+\mu_Y)\sqrt{x}}{1+\rho}\bigg)\big\{1+O(x^{-1/2})\big\}.
\end{align}
Using that $\mathrm{e}^{x}/2=\cosh(x)-\mathrm{e}^{-x}/2=\cosh(x)+O(x^{-1/2})$ as $x\rightarrow\infty$,
%that $\cosh(x)=(\mathrm{e}^x+\mathrm{e}^{-x})/2=\mathrm{e}^x/2+O(x^{-1/2})$ as $x\rightarrow\infty$, 
we deduce the asymptotic approximation (\ref{for1}).

To complete the proof of the asymptotic approximation (\ref{for1}) it remains to show that $R_1(x)/f_1(x)=O(x^{-1})$ as $x\rightarrow\infty$ and $R_2(x)/f_2(x)=O(x^{-1/2})$ as $x\rightarrow\infty$. We first argue that $R_2(x)/f_2(x)=O(x^{-1/2})$ as $x\rightarrow\infty$. Recall that in deriving the asymyptotic approximation (\ref{f2}) for $f_2(x)$ we used the limiting form (\ref{pfq}) to obtain that, as $x\rightarrow\infty$,
\begin{equation}\label{term1}
\sum_{n=0}^\infty\frac{n^{2k}}{(2n)!}y^n=\frac{1}{2^{2k+1}}y^k \mathrm{e}^{\sqrt{y}}\big\{1+O(x^{-1/2})\big\},  
\end{equation}
where $y=(\mu_X+\mu_Y)^2x/(1+\rho)^2$. Recall also the constants $c_k(n)$, $k,n\geq0$, that were introduced earlier in the proof. Since $c_0(n)=1$ and $c_k(n)$ is a polynomial in $n$ of degree $2k-1$ for $k\geq1$, instead of picking up the sum $\sum_{n=0}^\infty n^{2k}y^n/(2n)!$ we pick up a sum of the form $\sum_{n=0}^\infty C_k(n)y^n/(2n)!$, where $C_0(n)=1$ and, for $k\geq1$, $C_k(n)$ is a polynomial in $n$ of degree $2k-1$ with coefficients involving $k$ and $(\mu_X,\mu_Y,\rho)$. Applying the limiting form (\ref{pfq}) it follows that, as $x\rightarrow\infty$,
\begin{equation}\label{term2}
\sum_{n=0}^\infty\frac{C_k(n)}{(2n)!}y^n=\frac{\gamma_k}{2^{2k}}y^{k-1/2}\mathrm{e}^{\sqrt{y}}\big\{1+O(x^{-1/2})\big\},   
\end{equation}
where $\gamma_k$ (which depends on $k$ and $(\mu_X,\mu_Y,\rho)$ but not on $x$ and $n$) is the leading coefficient of the polynomial $C_k(n)$. On comparing (\ref{term1}) and (\ref{term2}), it follows that $R_2(x)/f_2(x)=O(x^{-1/2})$ as $x\rightarrow\infty$. It is now straightforward to see that $R_1(x)/f_1(x)=O(x^{-1})$ as $x\rightarrow\infty$. Since $b_k(m-n)$ is a polynomial in $n$ of degree $2k-2$, applying Lemma \ref{lem1} to $R_1(x)$ in the same way we did to $f_1(x)$ means that, as $x\rightarrow\infty$, $R_1(x)\sim R_3(x)+R_4(x)$ where $R_4(x)/R_3(x)=O(x^{-1/2})$ as $x\rightarrow\infty$, and $R_3(x)$ is of a similar form to $f_2(x)$, but rather than having a factor $n^{2k}$ in the sum has a factor that is a polynomial in $n$ of degree $2k-2$. By a similar argument to the one we used to show that $R_2(x)/f_2(x)=O(x^{-1/2})$ as $x\rightarrow\infty$, we deduce that $R_3(x)/R_1(x)=O(x^{-1})$  as $x\rightarrow\infty$, from which it follows that $R_1(x)/f_1(x)=O(x^{-1})$  as $x\rightarrow\infty$. This completes the proof of the asymptotic approximation (\ref{for1}).

Finally, we consider the limiting behaviour as $x\rightarrow-\infty$. 
%It is readily seen that the asymptotic approximations (\ref{}) and (\ref{}) follow from replacing $(x,\rho,\mu_X,\mu_Y)$ by $(-x,-\rho,\mu_X,-\mu_Y)$ in the asymptotic approximations (\ref{}) and (\ref{}) (in the case $\sigma_X^2=\sigma_Y^2=1$). We
We observe that $Z=XY=_d-X'Y'$, where $(X', Y')$ is a bivariate normal random vector with mean vector $(\mu_X,-\mu_Y)$, variances $(\sigma_X^2,\sigma_Y^2)$ and correlation coefficient $-\rho$. The asymptotic approximation (\ref{for2}) thus follows from replacing $(x,\rho,\mu_Y)$ by $(-x,-\rho,-\mu_Y)$ in the asymptotic approximation (\ref{for1}).
%(in the case $\sigma_X^2=\sigma_Y^2=1$). 
\begin{comment} 
 Applying the asymptotic expansion (\ref{Ktendinfinity}) to the PDF (\ref{pdf}) gives that, as $x\rightarrow-\infty$,
\begin{align}f(x)&\sim N(x)\sum_{n=0}^\infty\sum_{m=0}^{2n}\frac{(-x)^{n}}{(2n)!(1-\rho^2)^{2n}}\binom{2n}{m}(\rho \mu_Y-\mu_X)^m(\mu_Y-\rho \mu_X)^{2n-m}\nonumber\\
\label{second}&\quad\times\sum_{k=0}^\infty a_k(m-n)\frac{(1-\rho^2)^k}{(-x)^k},
\end{align} 
where
\begin{equation*}N(x)=\frac{1}{\sqrt{2\pi (-x)}}\exp\bigg\{-\frac{\mu_X^2+\mu_Y^2-2\rho\mu_X\mu_Y}{2(1-\rho^2)}\bigg\}\exp\bigg(\frac{x}{1-\rho}\bigg).
\end{equation*} 
We observe that the asymptotic expansion (\ref{second}) can be obtained from the asymptotic expansion (\ref{first}) by replacing $(x,\rho,\mu_X,\mu_Y)$ by $(-x,-\rho,\mu_X,-\mu_Y)$. The asymptotic approximation (\ref{for2}) now follows from replacing $(x,\rho,\mu_X,\mu_Y)$ by $(-x,-\rho,\mu_X,-\mu_Y)$ in the asymptotic approximation (\ref{for1}) (in the case $\sigma_X^2=\sigma_Y^2=1$). 
\end{comment}
\hfill $\Box$

\vspace{2mm}

\noindent{\emph{Proof of Theorem \ref{corsec3}.}} To obtain (\ref{for3}), we use that $\bar{F}(x)=\int_x^\infty f(t)\,\mathrm{d}t$ and then apply the limiting form (\ref{for1}) followed by Lemma \ref{lem2}, treating the cases $\mu_X+\mu_Y=0$ and $\mu_X+\mu_Y\not=0$ separately, as we did in the proof of Theorem \ref{prop3}. When dealing with the case $\mu_X+\mu_Y=0$ we first use that $\mathrm{cosh}(x)=\mathrm{e}^x/2+\mathrm{e}^{-x}/2=\mathrm{e}^x/2+O(x^{-1/2})$ as $x\rightarrow\infty$ in order to apply Lemma \ref{lem2}, and then after applying the lemma we use that $\mathrm{e}^x=2\cosh(x)-\mathrm{e}^{-x}=2\cosh(x)+O(x^{-1/2})$ as $x\rightarrow\infty$ to put the asymptotic approximation into the form (\ref{for3}). Note that the statements regarding the order of error in the asymptotic approximation (\ref{for3}) as $x\rightarrow\infty$ require part (ii) of Lemma \ref{lem2}. The proof of (\ref{for4}) is similar, this time making use of the formula $F(x)=\int_{-\infty}^x f(t)\,\mathrm{d}t$ and the limiting form (\ref{for2}), rather than (\ref{for1}). \hfill $\Box$
% Use that $\overline{F}(x)=\int_x^\infty f(t)\,\mathrm{d}t$ and $F(x)=\int_{-\infty}^x f(t)\,\mathrm{d}t$ and then 

\vspace{2mm}

\noindent{\emph{Proof of Theorem \ref{thmq}.}} We first derive the asymptotic approximation (\ref{q1}). The quantile function $Q(p)$ satisfies the equation $\bar{F}(Q(p))=1-p$. On applying the asymptotic approximation (\ref{for3}) we see that $Q(p)$ solves an equation of the form (\ref{zeqn}) from Lemma \ref{lemq} with
\begin{align*}
z&=1-p, \quad a=\frac{1}{\sigma_X\sigma_Y(1+\rho)},\quad b=\frac{1}{(1+\rho)\sqrt{\sigma_X\sigma_Y}}\bigg|\frac{\mu_X}{\sigma_X}+\frac{\mu_Y}{\sigma_Y}\bigg|,\\
A&=\frac{C(1+\rho)\sqrt{\sigma_X\sigma_Y}}{\epsilon\sqrt{2\pi}}\exp\bigg\{\frac{1}{8}\bigg(\frac{1+\rho}{1-\rho}\bigg)\bigg(\frac{\mu_X}{\sigma_X}-\frac{\mu_Y}{\sigma_Y}\bigg)^2\bigg\},
\end{align*}
where $\epsilon=1$ if $\mu_X/\sigma_X+\mu_Y/\sigma_Y=0$ and $\epsilon=2$ if $\mu_X/\sigma_X+\mu_Y/\sigma_Y\not=0$. Applying the asymptotic approximation (\ref{lemqeqn}) with these values of $z$, $a$, $b$ and $A$ and simplifying now yields the asymptotic approximation (\ref{q1}).

The derivation of the asymptotic approximation (\ref{q2}) is similar, this time using the fact that $Q(p)$ satisfies the equation $F(Q(p))=p$ and applying the asymptotic approximation (\ref{for4}); we omit the details. \hfill $\Box$ 

\vspace{2mm}

\noindent{\emph{Proof of Corollary \ref{cortvar}.}} Since $\mathrm{VaR}_p(Z)=F^{-1}(p)=Q(p)$, the asymptotic approximation (\ref{varfor}) is immediate from (\ref{q1}). Recall that $\mathrm{TVaR}_p(Z)=(1-p)^{-1}\int_p^1 \mathrm{VaR}_t(Z)\,\mathrm{d}t$. Applying this formula together with the asymptotic approximation (\ref{varfor}) and the formula $\int_p^1 \ln(1/(1-t))\,\mathrm{d}t=(1-p)[1+\ln(1/(1-p))]$
and the limiting forms
\begin{align*}
%\int_p^1 \ln(1/(1-t))\,\mathrm{d}t&=p\ln(1/p)+p,    \\
\int_p^1\sqrt{\ln(1/(1-t))}\,\mathrm{d}t&=(1-p)\sqrt{\ln(1/(1-p))} +O\bigg(\frac{1-p}{\sqrt{\ln(1/(1-p))}}\bigg), \quad p\rightarrow1, \\
\int_p^1\ln(\ln(1/(1-t)))\,\mathrm{d}t&=(1-p)\ln(\ln(1/(1-p)))+O\bigg(\frac{1-p}{\ln(1/(1-p))}\bigg), \quad p\rightarrow1, \\
\int_p^1\frac{\ln(\ln(1/(1-t)))}{\sqrt{\ln(1/(1-t))}}\,\mathrm{d}t&=O\bigg((1-p)\frac{\ln(\ln(1/(1-p)))}{\sqrt{\ln(1/(1-p))}}\bigg),\quad p\rightarrow1,\\
\int_p^1 \frac{1}{\ln(1/(1-t))}\,\mathrm{d}t&=O\bigg((1-p)\frac{1}{\ln(1/(1-p))}\bigg), \quad p\rightarrow1
\end{align*}
(which are easily obtained by integration by parts) yields (\ref{tvarfor}), as required. \hfill $\Box$

\section{Discussion}

Over the years, there has been much interest in the distribution of the product of correlated normal random variables with non-zero means, but it was not until the work of \cite{cui} in 2016 that a closed-form formula was obtained for the PDF.
This formula takes a rather complicated form. 
In this paper, we addressed the natural problem of deriving asymptotic approximations for the PDF. As a consequence, we were able to derive asymptotic approximations for tail probabilities, the quantile function, as well as asymptotic approximations for the risk measures VaR and TVaR. Our asymptotic approximations have yielded new theoretical insights into the asymptotic behaviour of these key distributional properties, and also allow for simple and computationally efficient approximations. 
%Our numerical results show that 

%Part of the motivation for studying the asymptotics of the PDF are that many important summary stat

Having described the leading order asymptotic behaviour of the PDF of the product $Z$, a natural follow on question is to ask for further terms in the asymptotic expansion; which could in turn be used to derive further terms in the asymptotic expansions of tail probabilities, the quantile function and TVaR. Further terms in the asymptotic expansions may lead to more accurate approximations that would still be simple and computationally efficient to implement. There is scope for extending the proof of Theorem \ref{prop3} to derive further terms in the asymptotic expansion, but the calculations would be rather involved, and this is left for a further work. 

We close by noting that very recently \cite{gnp24} obtained exact formulas for the PDF of the sum of independent copies of the product $Z=XY$ when the means of $X$ and $Y$ are non-zero; a closed-form formula for the PDF of $\overline{Z}_n$ in the zero-mean case $\mu_X=\mu_Y=0$ was previously known (see \cite{gaunt prod,man,np16}). The mean $\overline{Z}_n=n^{-1}(Z_1+Z_2+\cdots+Z_n)$, where $Z_1,Z_2,\ldots,Z_n$ are independent copies of $Z$, arises in areas such as astrophysics \cite{man,watts}, quantum cosmology \cite{g96} and electrical engineering \cite{ware}. %The PDF of $\overline{Z}_n$ in the zero-mean case $\mu_X=\mu_Y=0$ is known (see \cite{gaunt prod,man,np16}), and generalises the formula (\ref{pdf0}) in a simple manner with the modified Bessel function $K_0(|x|)$ replaced by $|x|^{(n-1)/2}K_{(n-1)/2}(|x|)$. It is therefore natural to conjecture that the PDF of $\overline{Z}_n$ can also be written as an infinite series involving the modified Bessel function of the second kind, 
It is possible that the methods used in this paper may prove useful in deriving asymptotic approximations for the PDF, tail probabilities, quantile function and TVaR of the mean $\overline{Z}_n$.
%We hope that our work provides motivation for researchers to address this open problem. 

%We close by noting that it is an open problem to find an exact closed-form formula for the PDF of the sum of independent copies of the product $Z=XY$ when the means of $X$ and $Y$ are non-zero. The mean $\overline{Z}_n=n^{-1}(Z_1+Z_2+\cdots+Z_n)$, where $Z_1,Z_2,\ldots,Z_n$ are independent copies of $Z$, arises in areas such as astrophysics \cite{man,watts}, quantum cosmology \cite{g96} and electrical engineering \cite{ware}. The PDF of $\overline{Z}_n$ in the zero-mean case $\mu_X=\mu_Y=0$ is known (see \cite{gaunt prod,man,np16}), and generalises the formula (\ref{pdf0}) in a simple manner with the modified Bessel function $K_0(|x|)$ replaced by $|x|^{(n-1)/2}K_{(n-1)/2}(|x|)$. It is therefore natural to conjecture that the PDF of $\overline{Z}_n$ can also be written as an infinite series involving the modified Bessel function of the second kind, in which case the methods used in this paper may prove useful in deriving asymptotic approximations for the PDF, tail probabilities, quantile function and TVaR. We hope that our work provides motivation for researchers to address this open problem.

\appendix

\section{Special functions}\label{appa}
In this appendix, we define the modified Bessel function of the second kind, 
%the modified Lommel function of the first kind, 
the complementary error function and the generalized hypergeometric function, and state some basic properties that are needed in this paper. 
%Unless  otherwise stated, t
These and further properties can be found in 
%the standard reference 
\cite{olver}. 

%\subsection{Modified Bessel function of the second kind}

The \emph{modified Bessel function of the second kind} is defined, for $\nu\in\mathbb{R}$ and $x>0$, by
\[K_\nu(x)=\int_0^\infty \mathrm{e}^{-x\cosh(t)}\cosh(\nu t)\,\mathrm{d}t.
\]
The modified Bessel function $K_\nu(x)$ has the following asymptotic behaviour:
\begin{eqnarray}\label{Ktend0}K_{\nu} (x) &\sim& \begin{cases} 2^{|\nu| -1} \Gamma (|\nu|) x^{-|\nu|}, & \quad x \downarrow 0, \: \nu \not= 0, \\
-\ln x, & \quad x \downarrow 0, \: \nu = 0, \end{cases} \\
\label{Ktendinfinity} K_{\nu} (x) &\sim& \sqrt{\frac{\pi}{2x}} \mathrm{e}^{-x}\sum_{k=0}^\infty\frac{a_k(\nu)}{x^k}, \quad x \rightarrow \infty,\: \nu\in\mathbb{R},
\end{eqnarray}
where $a_0(\nu)=1$ and
\begin{equation}\label{akv}a_k(\nu)=\frac{(4\nu^2-1^2)(4\nu^2-3^2)\cdots(4\nu^2-(2k-1)^2)}{k!8^k}, \quad k\geq1.
\end{equation}
\begin{comment}
We have the following indefinite integral formula
\begin{align*}\int x^\mu K_\nu(x)\,\mathrm{d}x=-x\big((\mu+\nu-1)K_\nu(x)t_{\mu-1,\nu-1}(x)+K_{\nu-1}t_{\mu,\nu}(x)\big),
\end{align*}
where $t_{\mu,\nu}(x)$ is a modified Lommel function of the first kind, defined by the hypergeometric series
\begin{align}t_{\mu,\nu}(x)&=\frac{x^{\mu+1}}{(\mu-\nu+1)(\mu+\nu+1)} {}_1F_2\bigg(1;\frac{\mu-\nu+3}{2},\frac{\mu+\nu+3}{2};\frac{x^2}{4}\bigg)\nonumber \\
&=2^{\mu-1}\Gamma\big(\tfrac{\mu-\nu+1}{2}\big)\Gamma\big(\tfrac{\mu+\nu+1}{2}\big)\sum_{k=0}^\infty\frac{(\frac{1}{2}x)^{\mu+2k+1}}{\Gamma\big(k+\frac{\mu-\nu+3}{2}\big)\Gamma\big(k+\frac{\mu+\nu+3}{2}\big)}. \nonumber
\end{align}
\end{comment}

The \emph{complementary error} function is defined, for $x\in\mathbb{R}$, by the integral
\begin{equation}\label{erfcdefn}\mathrm{erfc}(x)=\frac{2}{\sqrt{\pi}}\int_x^\infty\mathrm{e}^{-t^2}\,\mathrm{d}t.
\end{equation}
It has the following asymptotic behaviour as $x\rightarrow\infty$:
\begin{equation}\label{erfc}\mathrm{erfc}(x)=\frac{\mathrm{e}^{-x^2}}{\sqrt{\pi}x}\big(1+O(x^{-2})\big).
\end{equation}

The \emph{generalized hypergeometric function} is defined by the power series
\begin{equation}
\label{gauss}
{}_pF_q\bigg({a_1,\ldots,a_p \atop b_1,\ldots,b_q};x\bigg)=\sum_{j=0}^\infty\frac{(a_1)_j\cdots(a_p)_j}{(b_1)_j\cdots(b_q)_j}\frac{x^j}{j!},
\end{equation}
for $|x|<1$, and by analytic continuation elsewhere. Here $(u)_j=u(u+1)\cdots(u+k-1)$ is the ascending factorial. 

Let $k\geq1$. It is readily seen from the asymptotic expansion (16.11.8) of \cite{olver} that
\begin{align}\label{gausslim}{}_{2k-1}F_{2k}\bigg({2,\ldots,2 \atop 1,\ldots,1,\frac{3}{2}};\frac{x}{4}\bigg)= \frac{1}{2^{2k}}x^{k-1}\mathrm{e}^{\sqrt{x}}\big(1+O(x^{-1/2})\big), \quad x\rightarrow\infty.
\end{align}

\section*{Acknowledgements}
We would like to thank the reviewer for their helpful comments. RG is funded in part by EPSRC grant EP/Y008650/1 and EPSRC grant UKRI068. ZY is supported by a University of Manchester Research Scholar Award.

\end{document}